%% file: arxiv_version.tex
\documentclass[10pt,a4paper]{book}
\usepackage[utf8]{inputenc}
\usepackage{amsmath}
\usepackage{amsfonts}
\usepackage{amssymb}
\usepackage{authblk}
\usepackage{graphicx}

\usepackage{csquotes}

\usepackage[style=authoryear-ibid,backend=bibtex, maxcitenames = 2, mincitenames = 1, uniquelist=false, uniquename=false, abbreviate = true, maxbibnames = 10, minbibnames = 10]{biblatex}

\input{macro}

\title{A Tutorial on Bayesian Data Assimilation}

\author[1,2]{Colin Grudzien}
\author[3]{Marc Bocquet}
\affil[1]{Center for Western Weather and Water Extremes (CW3E), Scripps Institution of Oceanography, University of California San Diego, San Diego, CA, USA}
\affil[2]{Department of Mathematics and Statistics, University of Nevada, Reno, Reno, NV, USA}
\affil[3]{CEREA, \'{E}cole des Ponts and EDF R\&D, \^Ile-de-France, France}
\date{}
\begin{document}
\setcounter{chapter}{3}
\maketitle

\section*{Abstract}

This tutorial provides a broad introduction to Bayesian data assimilation that will be useful to practitioners, in interpreting algorithms and results, and for theoretical studies  developing novel schemes with an understanding of the rich history of geophysical data assimilation and its current directions.  The simple case of data assimilation in a `perfect' model is primarily discussed for pedagogical purposes.  Some mathematical results are derived at a high-level in order to illustrate key ideas about different estimators.  However, the focus of this work is on the intuition behind these methods, where more formal and detailed treatments of the data assimilation problem can be found in the various references.  In surveying a variety of widely used data assimilation schemes, the key message of this tutorial is how the Bayesian analysis provides a consistent framework for the estimation problem and how this allows one to formulate its solution in a variety of ways to exploit the operational challenges in the geosciences.

\section{Introduction}

In applications such as short- to medium-range weather prediction, \emph{data assimilation} (DA) provides a means to sequentially and recursively update forecasts of a time-varying physical process with newly incoming information \parencite{daley1991,kalnay2003,asch2016}, typically Earth observations.  The Bayesian approach to DA is widely adopted \parencite{lorenc1986} because it provides a unified treatment of tools from statistical estimation, nonlinear optimisation and machine learning for handling such a problem.  This tutorial illustrates how this approach can be utilised to develop and interpret a variety of widely used DA algorithms, both classical and those at the current state-of-the-art.

Suppose that the time-dependent physical states to be modelled can be written as a vector, $\pmb{x}_k \in \mathbb{R}^{N_x}$, where $k$ labels some time $t_k$. Formally, the time-evolution of these states is represented with the nonlinear map $\mathcal{M}$,
\begin{align}
\pmb{x}_k = \mathcal{M}_{k} \left(\pmb{x}_{k-1}, \boldsymbol{\lambda}\right) + \boldsymbol{\eta}_k
\label{eq:general_hmm}
\end{align}
where: (i) $\pmb{x}_{k-1}$ is the state variable vector at an earlier time $t_{k-1}$; (ii) $\boldsymbol{\lambda}$ is a vector of uncertain static physical parameters but on which the time evolution depends; and (iii) $\boldsymbol{\eta}_k$ is an additive (for simplicity but other choices are possible), stochastic noise term, representing errors in our model for the physical process.  Define $\Delta_t := t_{k} - t_{k-1}$ to be a fixed-length forecast horizon in this tutorial, though none of the following results require this to be fixed in practice. 
  
The basic goal of sequential DA is to estimate the random state vector $\pmb{x}_k$, given a prior distribution on $\left(\pmb{x}_{k-1}, \boldsymbol{\lambda}\right)$, and given knowledge of $\mathcal{M}_{k}$ and knowledge of how $\boldsymbol{\eta}_k$ is statistically distributed.  At time $t_{k-1}$, a forecast is made for the distribution of $\pmb{x}_k$ utilising the prior knowledge, which includes the physics-based model.   For simplicity, most of this tutorial is restricted to the case where $\boldsymbol{\lambda}$ is a known constant, and the forecast model is \emph{perfect}, i.e.,
\begin{align}
\pmb{x}_k = \mathcal{M}_{k} \left(\pmb{x}_{k-1}\right).
\end{align}
However, a strength of Bayesian analysis is how it easily extends to include a general treatment of model errors and the estimation of model parameters -- see, e.g., \textcite{asch2016} for a more general introduction.

While a forecast is made, one acquires a collection of observations of the real-world process.  This is written as the observation vector $\pmb{y}_k\in\mathbb{R}^{N_y}$, which is related to the state vector by
\begin{align}
\pmb{y}_k = \mathcal{H}_k \left(\pmb{x}_k\right) + \boldsymbol{\epsilon}_k.
\label{eq:obs_model}
\end{align}
The (possibly nonlinear) map $\mathcal{H}_k:\mathbb{R}^{N_x} \rightarrow \mathbb{R}^{N_y}$ relates the physical states being modelled, $\pmb{x}_k$, to the values that are actually observed, $\pmb{y}_k$.
Typically, in geophysical applications, observations are not $1:1$ with the state variables; while the data dimension can be extremely large, $N_y \ll N_x$ so this information is sparse relative to the model state dimension.  The term $\boldsymbol{\epsilon}_k$ in Eq.~\eqref{eq:obs_model} is an additive, stochastic noise term representing errors in the measurements, or a mismatch between the state variable representation and the observation \parencite[the representation error,][]{janjic2018}.

Therefore, some time after the real-life physical system has reached time $t_k$, one has a forecast distribution for the states $\pmb{x}_k$, generated by the prior on $\pmb{x}_{k-1}$ and the physics-based model $\mathcal{M}$, and the observations $\pmb{y}_k$ with some associated uncertainty.  The goal of Bayesian DA is to estimate the posterior distribution for $\pmb{x}_k$ conditioned on $\pmb{y}_k$, or some statistics of this distribution.
 
\section{Hidden Markov Models and Bayesian Analysis}
\subsection{The Observation-Analysis-Forecast Cycle}
Recursive estimation of the distribution for $\pmb{x}_k$ conditional on $\pmb{y}_k$ (i.e., assuming $\pmb{y}_k$ is known) can be described as an observation-analysis-forecast cycle \parencite{trevisan2004}.  Given the forecast-prior for the model state, and the likelihood function for the observation, Bayes' law updates the prior for the modelled state to the posterior conditioned on the observation.  Bayes' law is a simple re-arrangement of \emph{conditional probability}, defined by \textcite{kolmogorov2018foundations} as 
\begin{align}
\mathcal{P}(A|B):=\frac{\mathcal{P}(A,B)}{\mathcal{P}(B)},
\end{align}
for two events $A$ and $B$ where the probability of $B$ is non-zero.  Intuitively, this says that the probability of an event $A$, given knowledge that an event $B$ occurs, is equal to the probability of $A$ occurring relative to a sample space restricted to the event $B$.  Using the symmetry in the joint event $\mathcal{P}(A,B)=\mathcal{P}(B,A)$, Bayes' law is written
\begin{align}
\mathcal{P}(A|B) = \frac{\mathcal{P}(B|A) \mathcal{P}(A)}{\mathcal{P}(B)}.
\end{align}
In the observation-analysis-forecast cycle, $A$ is identified with the state vector (seen as a random vector) taking its value in a neighbourhood of $\pmb{x}_k$, $B$ is identified with the observation vector (seen as a random vector) taking its value in a neighbourhood of $\pmb{y}_k$.
The power of this statement is in how it describes an `inverse' probability -- while the \emph{posterior}, $\mathcal{P}(A|B)$,  on the left-hand-side may not be directly accessible, often the \emph{likelihood} $\mathcal{P}(B|A)$ and the \emph{prior} $\mathcal{P}(A)$ are easy to compute, and this is sufficient to develop a variety of probabilistic DA techniques.

\begin{figure}[ht]
\includegraphics[width=\linewidth]{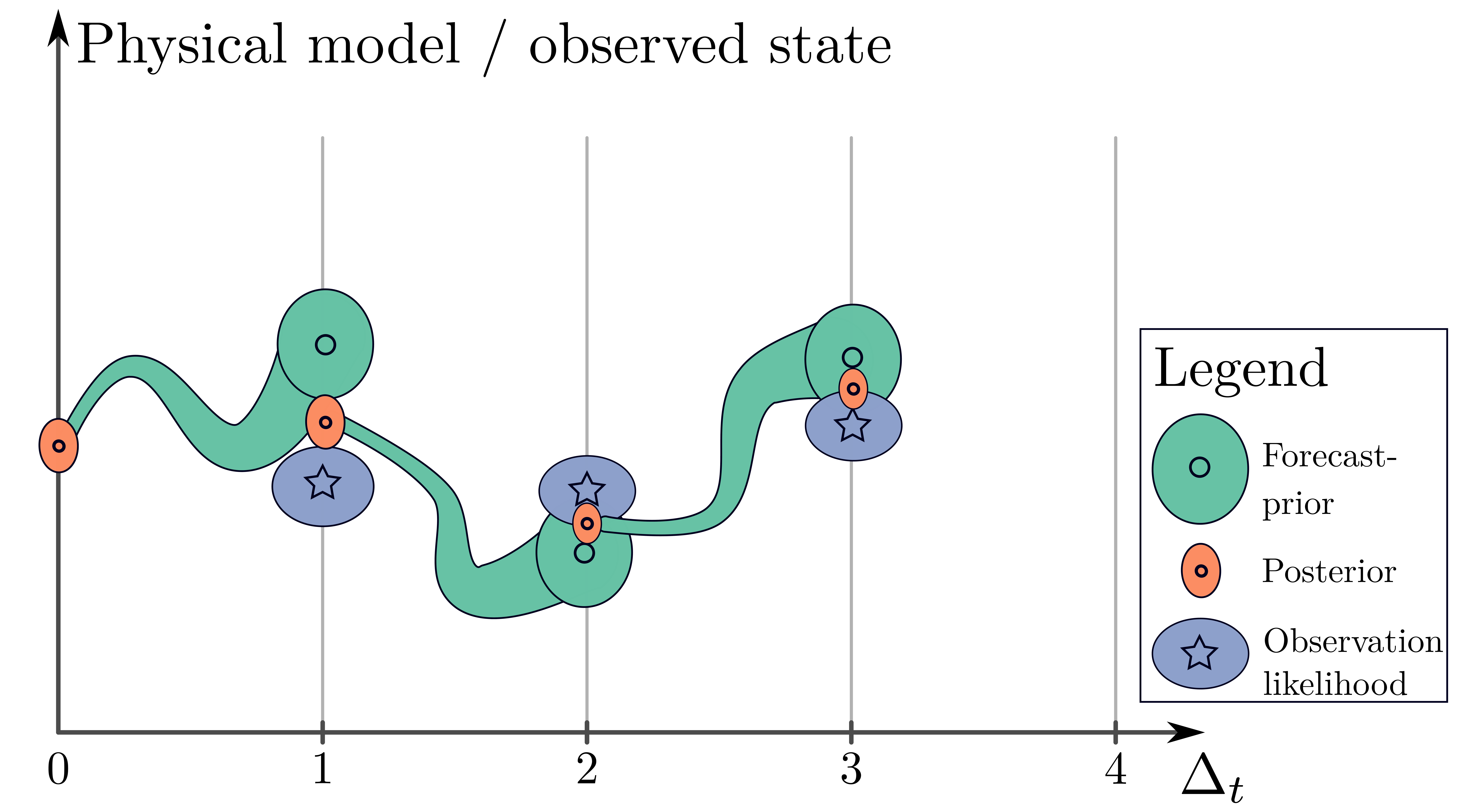}
\caption{Conceptual diagram of the observation-analysis-forecast cycle. The y-axis represents a variable of the state and observation vectors; the x-axis represents time.  Ellipses represent the spread of the observation / posterior / forecast errors.  Original figure, adapted from \textcite{carrassi2018}.}
\label{fig:filtering_diagram}
\end{figure}

A conceptual diagram of this process is pictured in Fig.~\ref{fig:filtering_diagram}.  Given the initial first prior (represented as the first `posterior' at time $t_0$), a forecast-prior for the model state at $t_1$,  $\mathcal{P}(A)$, is produced with the numerical model.  At the update time, there is an (possibly indirect and noisy) observation of the physical state with an associated likelihood $\mathcal{P}(B|A)$.  The posterior for the model state conditioned on this observation $\mathcal{P}(A|B)$, commonly denoted the \emph{analysis} in geophysical DA, is used to initialise the subsequent numerical forecast.  Recursive estimates of the current modelled state can be performed in this fashion, but a related question regards the past states.  The newly received observation gives information about the model states at past times and this allows one to produce a retrospective posterior estimate for the past states.  Recursive estimation of the present state using these incoming observations is commonly known as \emph{filtering}.  Conditional estimation of a past state given a time series including future observations is commonly known as \emph{smoothing}.

It is important to recognise that the filtering \emph{probability density function} (pdf) for the current time $p(\pmb{x}_k|\pmb{y}_k, \pmb{y}_{k-1}, \pmb{y}_{k-2}, \ldots)$ is actually just a marginal of the joint posterior pdf over all states in the current \emph{data assimilation window} (DAW), i.e., the window of lagged past and current states being estimated. In Fig.~\ref{fig:filtering_diagram}, the DAW is the time window $\{t_1,t_2,t_3\}$.  The conditional pdf for the model state at time $t_3$, given observations in the DAW is written in terms of the joint posterior as 
\begin{align}
\label{eq:filtering_marginal}
 p(\pmb{x}_3 \vert \pmb{y}_3, \pmb{y}_2, \pmb{y}_1) = \iiint \! p(\pmb{x}_3, \pmb{x}_2, \pmb{x}_1, \pmb{x}_0 \vert \pmb{y}_3, \pmb{y}_2, \pmb{y}_1) \, \mathrm{d}\pmb{x}_2 \mathrm{d}\pmb{x}_1 \mathrm{d}\pmb{x}_0,
\end{align}
by averaging out the past history of the model state from the joint posterior in the integrand.  A smoothing estimate may be produced in a variety of ways, exploiting different formulations of the Bayesian problem.  One may estimate only a \emph{marginal pdf} as on the left-hand-side of Eq.~\eqref{eq:filtering_marginal}, or the entire \emph{joint posterior pdf} as in the integrand above \parencite{anderson1979,cohn1994,cosme2012}.  The approach chosen may strongly depend on whether the DAW is static or is advanced in time.  Particularly, if one produces a smoothing estimate for the state at times $t_1$ through $t_3$, one may subsequently shift the DAW so that in the new cycle, the posterior for times $t_2$ through $t_4$ is estimated -- this type of analysis is known as \emph{fixed-lag smoothing}.  This tutorial considers how one can utilise a Bayesian \emph{maximum a posteriori} (MAP) formalism to efficiently solve the filtering and smoothing problems, using the various tools of statistical estimation, nonlinear optimisation and machine learning.

\subsection{A Generic Hidden Markov Model}

Recall the perfect physical process model, and the noisy observation model,
\begin{subequations}
\label{eq:perfect_model}
\begin{align}
\pmb{x}_k &= \mathcal{M}_{k} \left(\pmb{x}_{k-1}\right) \label{eq:state_model},\\
\pmb{y}_k &= \mathcal{H}_k \left(\pmb{x}_k\right) + \boldsymbol{\epsilon}_k. \label{eq:observation_model}
\end{align}
\end{subequations}
Denote the sequence of the process model states and observation model states between time $t_k$ and time $t_l$, for $k < l$, as
\begin{align}
\pmb{x}_{l:k} := \left\{\pmb{x}_l, \pmb{x}_{l-1}, \cdots, \pmb{x}_k\right\}, & &
\pmb{y}_{l:k} := \left\{\pmb{y}_l, \pmb{y}_{l-1}, \cdots, \pmb{y}_k\right\}. 
\end{align}
For arbitrary $l \in \mathbb{N}$, assume that the sequence of observation error 
\begin{align}
\{\boldsymbol{\epsilon}_l , \boldsymbol{\epsilon}_{l-1},\cdots, 1\}
\end{align}
is independent-in-time (i.e., a white process). 
  
The above formulation is a type of \emph{hidden Markov model}, where the dynamic state variables $\pmb{x}_k$ are known as the hidden variables because they are not directly observed.  A Markov model is a type of `memoryless' process, described this way because of how the conditional probability for the state is represented between different times \parencite[see Chapter 4]{ross2014introduction}. Particularly, if $\pmb{x}_{k:1}$ is a Markov process, the Markov property is defined as
\begin{align}
p\left(\pmb{x}_{k}|\pmb{x}_{k-1:0}\right) &= p\left(\pmb{x}_{k} \vert \pmb{x}_{k-1}\right) .
\label{eq:markov_assumption}
\end{align}
The above statement says that, given knowledge of the state $\pmb{x}_{k-1}$, the conditional probability for $\pmb{x}_k$ is independent of the past history of the state before time $t_{k-1}$, representing the probabilistic analogue of an initial value problem.  

Applying the Markov property recursively with the definition of the conditional pdf yields
\begin{align}
p\left(\pmb{x}_{L:0}\right) = p(\pmb{x}_0) \prod_{k=1}^{L} p(\pmb{x}_k|\pmb{x}_{k-1}).
\end{align}
Therefore, the joint pdf for a forecast of the model state can be written as the product of the first prior at time $t_0$, representing the uncertainty of the data used to initialise the model forecast, and the product of the \emph{Markov transition pdfs}, describing the evolution of the state between discrete times.

With the perfect state model, as in Eq.~\eqref{eq:state_model} above, the transition probability for some subset $\mathrm{d}\pmb{x}\subset \mathbb{R}^{N_x}$ is written
\begin{align}
\mathcal{P}\left(\pmb{x}_k \in \mathrm{d}\pmb{x} \vert \pmb{x}_{k-1}\right) = \delta_{\mathcal{M}_k\left(\pmb{x}_{k-1}\right)}(\mathrm{d}\pmb{x}),
\end{align} 
with $\delta_{\pmb{v}}$ referring to the \emph{Dirac measure} at $\pmb{v} \in \mathbb{R}^{N_x}$. The Dirac measure satisfies
\begin{align}
\int \! f(\pmb{x}) \delta_{\pmb{v}}\left(\mathrm{d}\pmb{x}\right) = f\left(\pmb{v}\right),
\end{align}
where this is a singular measure, to be understood by the integral equation.  Accordingly, the transition pdf is often written proportional as
\begin{align}
  p(\pmb{x}_k \vert \pmb{x}_{k-1}  ) \propto \delta \left\{\pmb{x}_k - \mathcal{M}_k\left(\pmb{x}_{k-1}\right)\right\}
\end{align}
where $\delta$ represents the \emph{Dirac distribution}. Heuristically, this is known as the `function' defined by the property
\begin{align}
 \delta_\epsilon(\pmb{x}) = \begin{cases} \frac{1}{\epsilon} & \pmb{x} \in \left[-\epsilon,+ \epsilon\right] \\
    0 & \text{else}\end{cases},
    & & \delta(\pmb{x}) = \lim_{\epsilon \rightarrow 0^+} \delta_{\epsilon}(\pmb{x}).
\end{align}
However, this is just a convenient abuse of notations, as the Dirac measure does not have a pdf with respect to the standard Lebesgue measure. Rather, the Dirac distribution is understood through the generalised function theory of distributions \parencite[see Section 3.4]{taylor1996partial} as a type of integral kernel satisfying
\begin{align}
  \int \! f(\pmb{x}_{k}) \delta\left\{\pmb{x}_k - \mathcal{M}_k\left(\pmb{x}_{k-1}\right)\right\}\mathrm{d}\pmb{x}_{k} = f\left(\mathcal{M}_k\left(\pmb{x}_{k-1}\right)\right).\label{eq:dirac_distribution}
\end{align}

Equation~\eqref{eq:dirac_distribution} represents the existence and uniqueness of the solution to an initial value problem in deterministic systems of ordinary and partial differential equations, where the knowledge of the state $\pmb{x}_{k-1}$ completely determines the state $\pmb{x}_k$ via the perfect forecast model $\mathcal{M}_k$.  Particularly, there is probability $1$ of the state following the unique solution to the time evolution, and probability $0$ of all other outcomes. Each Markov transition pdf represents the evolution of an initial condition with respect to the dynamical model, conditional ultimately on an uncertain outcome of the initial data from the first prior.  While the perfect model assumption is used for simplicity, the decomposition of the forecast pdf above can be derived for erroneous models under the additional assumption that the model errors are independent-in-time.  Like the decomposition of the forecast pdf,  for the observation likelihood we can write
\begin{align}
p\left(\pmb{y}_{k}\vert \pmb{x}_k, \pmb{y}_{k-1:1}\right) = p\left(\pmb{y}_k \vert \pmb{x}_k\right)
\label{eq:observation_independence}
\end{align}
due to the independence assumption on the observation errors, and the relationship between $\pmb{x}_k$ and $\pmb{y}_k$.  This says that the knowledge of the physical state $\pmb{x}_k$ completely determines the likelihood of the observation $\pmb{y}_k$.  Indeed,
\begin{align}
\pmb{\epsilon}_k = \pmb{y}_k - \mathcal{H}_k(\pmb{x}_k),
\end{align}
which follows a known distribution and is independent of the other observation error outcomes by assumption.

Consider thus how to estimate the filtering pdf $p\left(\pmb{x}_{k} \vert \pmb{y}_{k:1}\right)$. Using the definition of the conditional pdf, one has
\begin{align}
 p\left(\pmb{x}_{k} \vert \pmb{y}_{k:1}\right) &= \frac{p\left(\pmb{y}_{k:1},\pmb{x}_k \right)}{ p\left(\pmb{y}_{k:1}\right)}.
\end{align}
Rewriting these pdfs as conditional pdfs, and by using the independence assumption,
\begin{subequations}
\begin{align}
p\left(\pmb{x}_{k} \vert \pmb{y}_{k:1}\right) &=\frac{p\big(\pmb{y}_k, (\pmb{x}_k, \pmb{y}_{k-1:1})\big)}{p\left(\pmb{y}_{k:1}\right)}\\
&=\frac{p\left(\pmb{y}_{k}\vert \pmb{x}_k, \pmb{y}_{k-1:1}\right) p\left(\pmb{x}_k, \pmb{y}_{k-1:1}\right)}{p\left(\pmb{y}_{k:1}\right)}  =\frac{p\left(\pmb{y}_{k}\vert \pmb{x}_k\right) p\left(\pmb{x}_k, \pmb{y}_{k-1:1}\right)}{p\left(\pmb{y}_{k:1}\right)}.
\end{align}
\end{subequations}
Writing the joint pdfs again in terms of conditional pdfs
\begin{subequations}
\begin{align}
p\left(\pmb{x}_{k} \vert \pmb{y}_{k:1}\right) &=\frac{p\left(\pmb{y}_k \vert \pmb{x}_k\right) p\left(\pmb{x}_k\vert \pmb{y}_{k-1:1}\right) p\left(\pmb{y}_{k-1:1}\right)}{p\left(\pmb{y}_k\vert \pmb{y}_{k-1:1}\right) p\left(\pmb{y}_{k-1:1}\right)} \\ 
&=\frac{p\left(\pmb{y}_k \vert \pmb{x}_k\right) p\left(\pmb{x}_k\vert   \pmb{y}_{k-1:1}\right)}{p\left(\pmb{y}_k\vert \pmb{y}_{k-1:1}\right)}.
\end{align}
\label{eq:filtering_recursion}
\end{subequations}

Now, suppose that the posterior pdf $p(\pmb{x}_{k-1} \vert \pmb{y}_{k-1:1})$ at the last observation time $t_{k-1}$ is already computed -- then the model forecast of this pdf is given by averaging over the state at time $t_{k-1}$ with respect to the Markov transition pdf,
\begin{align}
p(\pmb{x}_k\vert \pmb{y}_{k-1:1}) &= \int \! p(\pmb{x}_k \vert \pmb{x}_{k-1}) p(\pmb{x}_{k-1} \vert \pmb{y}_{k-1:1}) \, \mathrm{d}\pmb{x}_{k-1},
\end{align}
yielding the forecast-prior.  The filtering pdf, on the left-hand-side of Eq.~\eqref{eq:filtering_recursion} is written in terms of: (i) the likelihood of the observed data given the model forecast, $p\left(\pmb{y}_{k}\vert \pmb{x}_{k}\right)$; (ii) the forecast-prior given the last best estimate of the state, $p\left(\pmb{x}_k\vert   \pmb{y}_{k-1:1}\right)$; and (iii) the marginal of the joint pdf $p( \pmb{y}_k, \pmb{x}_k \vert \pmb{y}_{k-1:1})$, integrating out the hidden variables,
\begin{align}
p\left(\pmb{y}_k \vert \pmb{y}_{k-1:1}\right) = \int \! p(\pmb{y}_k \vert \pmb{x}_k) p(\pmb{x}_k \vert \pmb{y}_{k-1:1}) \, \mathrm{d}\pmb{x}_{k}.
\end{align}
This type of pdf, only depending on the observations, is called an \emph{evidence} \parencite[e.g.,][]{carrassi2017}.

Typically, the pdf in the denominator of Eq.~\eqref{eq:filtering_recursion} is mathematically intractable. However, the denominator is independent of the hidden variable $\pmb{x}_{k}$ by construction -- the free argument in the pdf on the left-hand-side is the model state $\pmb{x}_k$ and the purpose of the denominator on the right-hand-side is only to normalise the integral of the posterior pdf to $1$.  Instead, as a proportionality statement,
\begin{align}
p\left(\pmb{x}_{k} \vert \pmb{y}_{k:1}\right) &\propto  p\left(\pmb{y}_k \vert \pmb{x}_k\right)  p\left(\pmb{x}_k\vert   \pmb{y}_{k-1:1}\right) 
\label{eq:filtering_propto}
\end{align}
one can devise the Bayesian MAP estimate as the choice of $\overline{\pmb{x}}_k$ that maximises the posterior pdf, but written in terms of the two right-hand-side components in Eq.~\eqref{eq:filtering_propto}.  For the purpose of maximising the posterior pdf, the denominator leads to insignificant constants that can be discarded.  Thus, in order to compute the MAP sequentially and recursively in time, one can develop a recursion in proportionality as above. 
However, note that the evidence can be estimated and used for other significant purposes still within a Bayesian framework \parencite{carrassi2017}.

\subsection{Linear-Gaussian Models}
 Generally, the filtering pdf $p\left(\pmb{x}_{k} \vert \pmb{y}_{k:1}\right) $ has no analytical solution, i.e., no explicit expression.  However, when the models are linear, i.e. both the state and observation models are written as matrix actions
\begin{subequations}
\begin{align}
\pmb{x}_k &= \mathbf{M}_k \pmb{x}_{k-1},\\
\pmb{y}_k &= \mathbf{H}_k \pmb{x}_k + \pmb{\epsilon}_k,
\end{align}
\label{eq:linear_gaussian}
\end{subequations}
and the error pdfs are \emph{Gaussian}:
\begin{align}
p(\pmb{x}_0) &= n\left(\pmb{x}_0 | \overline{\pmb{x}}_0, \mathbf{B}_0\right), \quad p\left(\pmb{y}_k \vert \pmb{x}_k \right) = n\left(\pmb{y}_k | \mathbf{H}_k \pmb{x}_k, \mathbf{R}_k\right),\\
n(\pmb{z}| \overline{\pmb{z}}, \mathbf{B}) &:= \frac{1}{ \sqrt{\left(2\pi\right)^{N_z} \det\left(\mathbf{B}\right)}}\exp\left\{ -\frac{1}{2}\left(\overline{\pmb{z}} - \pmb{z}\right)^\top\mathbf{B}^{-1}\left(\overline{\pmb{z}} - \pmb{z}\right) \right\},\label{eq:multivariate_gaussian_density}
\end{align}
then the forecast and posterior pdfs are Gaussian at all times, and are parametrised in terms of their mean and covariance.  In particular, Gaussian distributions are closed under affine transformations, i.e., maps of the form
\begin{align}
\pmb{f}(\pmb{x}) = \mathbf{A}\pmb{x} + \pmb{b},
\end{align}
corresponding to a linear transformation when $\pmb{b}$ is a vector of zeros, i.e., $\pmb{b}=\pmb{0}$ \parencite[see Theorem 3.3.3]{tong2012multivariate}.  If $\pmb{x}$ is distributed with pdf $p(\pmb{x} )= n\left(\pmb{x}| \overline{\pmb{x}}, \mathbf{B}\right)$, then the random vector $\pmb{y} :=  \mathbf{A}\pmb{x} + \pmb{b}$ is distributed with pdf
\begin{align}
p(\pmb{y}) = n \left(\pmb{y}|\mathbf{A}\overline{\pmb{x}} + \pmb{b}, \mathbf{A}\mathbf{B}\mathbf{A}^\top \right).
\end{align}
Suppose that the last analysis pdf is given as
\begin{align}
p(\pmb{x}_{k-1} | \pmb{y}_{k-1:1}) = n\left(\pmb{x}_{k-1} | \overline{\pmb{x}}_{k-1}^\mathrm{a}, \mathbf{B}_{k-1}^\mathrm{a}\right),
\end{align} 
parametrised in terms of the analysis mean $\overline{\pmb{x}}^\mathrm{a}_{k-1}$ and analysis error covariance $\mathbf{B}^\mathrm{a}_{k-1}$. Then, the forecast-prior pdf is written as
\begin{align}
p(\pmb{x}_k| \pmb{y}_{k-1:1}) = n\left(\pmb{x}_k |  \overline{\pmb{x}}^\mathrm{f}_{k-1},  \mathbf{B}^\mathrm{f}_{k} \right),
\end{align}
where the forecast mean and forecast error covariance are defined by
\begin{subequations}
\begin{align}
\overline{\pmb{x}}_k^\mathrm{f} := \mathbf{M}_k \overline{\pmb{x}}^\mathrm{a}_{k-1}, \\
\mathbf{B}_k^\mathrm{f}:= \mathbf{M}_k \mathbf{B}^\mathrm{a}_{k-1} \mathbf{M}_k^\top,
\end{align}
\label{eq:kf_forecast}
\end{subequations}
respectively.

Similarly, the conditional and marginal distributions of a Gaussian random vector are also Gaussian and their pdf has an analytical form. Suppose that $\pmb{z}\in \mathbb{R}^{N_z}$ is an arbitrary Gaussian random vector, partitioned as
\begin{align}
\pmb{z}:= \begin{pmatrix} \pmb{x} \\ \pmb{y} \end{pmatrix}  & & p(\pmb{z}):= n\left(\pmb{z} \Big|\begin{pmatrix}\overline{\pmb{x}} \\ \overline{\pmb{y}}\end{pmatrix}, 
\begin{pmatrix}
\boldsymbol{\Sigma}_{xx} & \boldsymbol{\Sigma}_{xy} \\
\boldsymbol{\Sigma}_{yx} & \boldsymbol{\Sigma}_{yy}
\end{pmatrix}\right),
\end{align}
with the dimensions given as $N_z = N_x + N_y$ and
\begin{subequations}
\begin{align}
\pmb{x},\overline{\pmb{x}} \in \mathbb{R}^{N_x}, & & \pmb{y},\overline{\pmb{y}} \in \mathbb{R}^{N_y}, & & \boldsymbol{\Sigma}_{xx} \in \mathbb{R}^{N_x\times N_x}, \\
\boldsymbol{\Sigma}_{xy}=\boldsymbol{\Sigma}_{yx}^\top \in \mathbb{R}^{N_x \times N_y},  & & \boldsymbol{\Sigma}_{yy} \in \mathbb{R}^{N_y \times N_y}. & &
\end{align}
\end{subequations}
Then, the general form of the pdf for $\pmb{x} $ conditioned on the outcome of $\pmb{y}$ is given by the Gaussian
\begin{align}
p(\pmb{x}|\pmb{y}) = n\left(\pmb{x}|\overline{\pmb{x}} + \boldsymbol{\Sigma}_{xy}\boldsymbol{\Sigma}_{yy}^{-1}\left(\pmb{y} -  \overline{\pmb{y}}\right), \boldsymbol{\Sigma}_{xx} - \boldsymbol{\Sigma}_{xy} \boldsymbol{\Sigma}_{yy}^{-1} \boldsymbol{\Sigma}_{yx}\right),
\end{align}
where the covariance matrix $\boldsymbol{\Sigma}_{xx} - \boldsymbol{\Sigma}_{xy} \boldsymbol{\Sigma}_{yy}^{-1} \boldsymbol{\Sigma}_{yx}$ is called the \emph{Schur complement} \parencite[see, e.g., Theorem 3.3.4 of][]{tong2012multivariate}.
Noting the form of the linear observation model in Eq.~\eqref{eq:linear_gaussian}, and the independence of the observation errors, it is easy to see that the vector composed of $\begin{pmatrix} \pmb{x}_k^\top & \pmb{y}_k^\top\end{pmatrix}^\top$ is jointly Gaussian.  Relying on the identifications
\begin{align}
\boldsymbol{\Sigma}_{xy} = \mathbf{B}_k^\mathrm{f}\mathbf{H}_k^\top, & & \boldsymbol{\Sigma}_{yy}= \mathbf{H}_k\mathbf{B}_k^\mathrm{f}\mathbf{H}_k^\top + \mathbf{R}_k,
\end{align}
the posterior pdf at time $t_k$ is derived as the Gaussian with analysis mean and analysis error covariance given by
\begin{subequations}
\label{eq:kalman-update}
\begin{align}
\overline{\pmb{x}}_k^\mathrm{a} &:= \overline{\pmb{x}}_k^\mathrm{f} + \mathbf{B}_k^\mathrm{f}\mathbf{H}_k^\top \left(\mathbf{H}_k\mathbf{B}_k^\mathrm{f}\mathbf{H}_k^\top + \mathbf{R}_k\right)^{-1}\left(\pmb{y}_k - \mathbf{H}_k\overline{\pmb{x}}_k^\mathrm{f}\right), \\
\mathbf{B}_k^\mathrm{a} &:= \mathbf{B}_k^\mathrm{f} - \mathbf{B}_k^\mathrm{f}\mathbf{H}_k^\top \left(\mathbf{H}_k\mathbf{B}_k^\mathrm{f}\mathbf{H}_k^\top + \mathbf{R}_k\right)^{-1} \mathbf{H}_k \mathbf{B}_k^\mathrm{f}.
\end{align}
\end{subequations}
Defining $\mathbf{K}_k:= \mathbf{B}_k^\mathrm{f}\mathbf{H}_k^\top \left(\mathbf{H}_k\mathbf{B}_k^\mathrm{f}\mathbf{H}_k^\top + \mathbf{R}_k\right)^{-1}$ as the Kalman gain, Eqs.~\eqref{eq:kalman-update} yield the classical \emph{Kalman filter} (KF) update, and this derivation inductively defines the forecast and posterior distribution for $\pmb{x}_k$ at all times.  The KF is recognised thus as the parametric representation of a linear-Gaussian hidden Markov model, providing a recursion on the first two moments for the forecast and posterior.  

This analysis extends, as with the classical KF, easily to incorporate additive model errors as in Eq.~\eqref{eq:general_hmm}, see e.g., \textcite{anderson1979}.  However, there are many ways to formulate the KF and there are some drawbacks of the above approach.  Even when the model forecast equations themselves are linear, if they functionally depend on an uncertain parameter vector, defined $\mathbf{M}_{k}(\pmb{\lambda})$, the joint estimation problem of $(\pmb{x}_k, \pmb{\lambda})$ can become highly nonlinear and an iterative approach to the joint estimation may be favourable.  Therefore, this tutorial develops the subsequent extensions to the KF with the MAP approach, which coincides with the development in least-squares and nonlinear optimisation.

\section{Least-Squares and Nonlinear Optimisation}
\subsection{The 3D Cost Function from Gaussian Statistics}

As seen in the previous section, the linear-Gaussian analysis admits an analytical solution for the forecast and posterior pdf at all times.  However, an alternative approach for deriving the KF is formed using the MAP estimation in proportionality.  Consider that the \emph{natural logarithm} ($\log$) is monotonic, so that an increase in the input of the argument corresponds identically to an increase in the output of the function.   Therefore, maximising the posterior pdf as in Eq.~\eqref{eq:filtering_propto} is equivalent to maximising
\begin{align}
\log\left( p\left(\pmb{y}_{k} \vert \pmb{x}_k\right)  p\left(\pmb{x}_k\vert   \pmb{y}_{k-1:1}\right) \right)& = \log\left( p\left(\pmb{y}_k \vert \pmb{x}_k\right) \right) + \log\left(  p\left(\pmb{x}_k|\pmb{y}_{k-1:1}\right)\right).
\end{align}
Given the form of the the multivariate Gaussian pdf in Eq.~\eqref{eq:multivariate_gaussian_density}, maximising the above log-posterior is equivalent to minimising the following least-squares cost function, derived in proportionality to the minus-log-posterior:
\begin{align}
\mathcal{J}_{\mathrm{KF}}(\pmb{x}_{k}) &= \frac{1}{2} \parallel \overline{\pmb{x}}_{k}^\mathrm{f} -\pmb{x}_{k} \parallel_{\mathbf{B}_{k}^\mathrm{f}}^2 + \frac{1}{2}\parallel \pmb{y}_k - \mathbf{H}_k \pmb{x}_k\parallel_{\mathbf{R}_k}^2.\label{eq:simple_map}
\end{align}
For an arbitrary positive definite matrix $\mathbf{A}$ the \emph{Mahalanobis distance} \parencite{mahalanobis1936generalized} with respect to $\mathbf{A}$ is defined as
\begin{align}
\parallel \pmb{v} \parallel_\mathbf{A} := \sqrt{\pmb{v}^\top \mathbf{A}^{-1}\pmb{v}}.
\end{align}
The above distances are weighted Euclidean norms with: (i) $\parallel \circ \parallel_{\mathbf{B}_k^\mathrm{f}}$ weighting relative to the  forecast spread; and (ii) $\parallel \circ \parallel_{\mathbf{R}_k}$ weighting relative to the observation imprecision.  The MAP state thus interpolates the forecast mean and the observation relative to the uncertainty in each piece of data.  Due to the unimodality of the Gaussian, and its symmetry about its mean, it is clear that the conditional mean $\overline{\pmb{x}}^\mathrm{a}_k$ is also the MAP state.

While this cost function analysis provides the solution to finding the first moment of the Gaussian posterior, this does not yet address how to find the posterior error covariance.  In this linear-Gaussian setting, it is easily shown that the analysis error covariance is actually given by
\begin{align}
\begin{split}
\mathbf{B}_k^\mathrm{a}:=&\left[ \left( \mathbf{B}_k^\mathrm{f} \right)^{-1} + \mathbf{H}_k^\top \mathbf{R}_k^{-1}\mathbf{H}_k \right]^{-1}
= \boldsymbol{\Xi}_{\mathcal{J}_\mathrm{KF}}^{-1},
\end{split}
\end{align}
where $\boldsymbol{\Xi}_{\mathcal{J}_\mathrm{KF}}$ refers to the \emph{Hessian} of the least-squares cost function in Eq.~\eqref{eq:simple_map}, i.e., the matrix of its mixed second partial derivatives in the state vector variables.  This is a fundamental result that links the recursive analysis in time to the MAP estimation performed with linear least-squares, see, e.g., Section 6.2 of \textcite{reich_cotter_2015} for further details.  

A more general result from maximum likelihood estimation extends this analysis as an approximation to nonlinear state and observation models, and to non-Gaussian error distributions.  Define $\hat{\pmb{\theta}}$ to be the \emph{maximum likelihood estimator} (MLE) for some unknown parameter vector $\pmb{\theta}\in \mathbb{R}^{N_\theta}$, where $\hat{\pmb{\theta}}$ depends on the realisation of some arbitrarily distributed random sample of observed data $\{\pmb{z}_i\}_{i=1}^N$. Then, under fairly general regularity conditions, the MLE satisfies
\begin{align}
\mathcal{I}\left(\pmb{\theta}\right)^{-\frac{1}{2}}\left(\hat{\pmb{\theta}} - \pmb{\theta}\right) \rightarrow_{d} \mathcal{N}\left(\pmb{0}, \mathbf{I}_{N_\theta}\right)
\label{eq:asymptotic_mle}
\end{align}
where: (i) Eq.~\eqref{eq:asymptotic_mle} refers to convergence in distribution of the random vector $\mathcal{I}\left(\pmb{\theta}\right)^{-\frac{1}{2}}\left(\hat{\pmb{\theta}} - \pmb{\theta}\right)$ as the sample size $N\rightarrow \infty$; (ii) $\mathbf{I}_{N_\theta}$ is the identity matrix in the dimension of $\pmb{\theta}$; (iii) $\mathcal{N}(\pmb{0},\mathbf{I}_{N_\theta})$ is the multivariate Gaussian distribution with mean zero and covariance equal to the identity matrix, i.e., with pdf $n(\pmb{x}|\pmb{0},\mathbf{I}_{N_\theta})$; and (iv) $\mathcal{I}(\pmb{\theta})$ is the \emph{Fisher information matrix}.  The Fisher information matrix is defined as the expected value of the Hessian of the minus-log-likelihood, taken over the realisations of the observed data $\{\pmb{z}_i\}_{i=1}^N$, and with respect to the true parameter vector $\pmb{\theta}$. It is common to approximate the distribution of the MLE as
\begin{align}
\hat{\pmb{\theta}} \sim \mathcal{N}\left(\pmb{\theta}, I\left(\pmb{\hat{\theta}}\right) \right)
\end{align}
where $I\left(\hat{\pmb{\theta}}\right)$ is the observed Fisher information, i.e., the realisation of the Hessian of the minus-log-likelihood given the observed sample data. This approximation improves in large sample sizes like the central limit theorem, see, e.g, Chapter 9 of \textcite{pawitan2001} for the details of the above discussion.  In nonlinear optimisation, the approximate distribution of the MLE therefore relates the geometry in the neighbourhood of of a local minimiser to the variation in the optimal estimate.  Particularly,  the curvature of the cost function level contours around the local minimum is described by the Hessian while the spread of the distribution of the MLE is described by the analysis error covariance.

\subsection{3D-VAR and the Extended Kalman Filter}
One of the benefits of the above cost function approach is that this immediately extends to handle a nonlinear observation operator $\mathcal{H}_k$ as a problem of nonlinear least-squares.  When the observations are related nonlinearly to the state vector, the Bayesian posterior is no longer generally Gaussian, but the interpretation of the analysis state $\overline{\pmb{x}}^\mathrm{a}_k$ interpolating  between a background proposal $\overline{\pmb{x}}_k^\mathrm{f}$ and the observed data $\pmb{y}_k$ relative to their respective uncertainties remains valid.  The use of the nonlinear least-squares cost function can be considered as making a Gaussian approximation, similar to the large sample theory in maximum likelihood estimation.  However, the forecast and analysis background error covariances, $\mathbf{B}^\mathrm{f/a}_k$, take on different interpretations depending on the DA scheme.

In full-scale geophysical models, the computation and storage of the background error covariance $\mathbf{B}_k^\mathrm{f/a}$ is rarely feasible due to its large size -- in practice, this is usually treated abstractly in its action as an operator by preconditioning the optimisation \parencite{tabeart2018preconditioning}.  One traditional approach to handle this reduced representation is to use a static-in-time background for the cost function in Eq.~\eqref{eq:simple_map}, rendering the \emph{three-dimensional, variational} (3D-VAR) cost function
\begin{align}
\mathcal{J}_{\mathrm{3D}\text{-}\mathrm{VAR}}(\pmb{x}) := \frac{1}{2}\parallel \overline{\pmb{x}}_k^\mathrm{f}  - \pmb{x} \parallel_{\mathbf{B}_{\mathrm{3D}\text{-}\mathrm{VAR}}}^2 + \frac{1}{2} \parallel \pmb{y}_k - \mathcal{H}_k\left(\pmb{x}\right)\parallel^2_{\mathbf{R}_k}.
\label{eq:3dvar}
\end{align}
The 3D-VAR background error covariance is typically defined as a posterior, `climatological' covariance, taken with respect to a long-time average over the modelled state.  One way that this is roughly estimated in meteorology is by averaging over a \emph{reanalysis data set}, which is equivalent to averaging with respect to a joint smoothing posterior over a long history for the system \parencite[and references therein]{kalnay2003}.  If one assumes that the DA cycle represents a stationary, ergodic process in its posterior statistics, and that this admits an invariant, ergodic measure, $\pi^\ast$, one would define the 3D-VAR background error covariance matrix with the expected value with respect to this measure as
\begin{subequations}
\begin{align}
\overline{\pmb{x}}^\ast &:= \mathbb{E}_{\pi^\ast}\left[\pmb{x} \right],\\
\mathbf{B}_{\mathrm{3D}\text{-}\mathrm{VAR}} &:= \mathbb{E}_{\pi^\ast}\left[ \left(\pmb{x} - \overline{\pmb{x}}^\ast \right) \left(\pmb{x} - \overline{\pmb{x}}^\ast \right)^\top\right].
\end{align}
\end{subequations}
The 3D-VAR cycle is defined where $\overline{\pmb{x}}^\mathrm{a}_k$ is the \emph{minimising argument (argmin)} of Eq.~\eqref{eq:3dvar} and subsequently $\overline{\pmb{x}}_{k+1}^\mathrm{f}:= \mathcal{M}_{k+1}\left(\overline{\pmb{x}}^\mathrm{a}_k\right)$, i.e., the forecast mean is estimated with the background control trajectory propagated through the fully nonlinear model.  Therefore, at every iteration of the algorithm, 3D-VAR can be understood to treat the proposal $\overline{\pmb{x}}^\mathrm{f}_k$ as a random draw from the invariant, climatological-posterior measure of the system; the optimised solution $\overline{\pmb{x}}_k^\mathrm{a}$ is thus the state that interpolates the forecast mean and the currently observed data as if the forecast mean was drawn randomly from the climatological-posterior.

For illustration, consider a cost-effective, explicit formulation of 3D-VAR where the background error covariance matrix is represented as a Cholesky factorisation.  Here, one performs the optimisation by writing the modelled state as a perturbation of the forecast mean,
\begin{align}
\pmb{x}_k := \overline{\pmb{x}}_k^\mathrm{f} + \boldsymbol{\Sigma} \pmb{w},
\label{eq:weight_definition}
\end{align}
where $\mathbf{B}_{\mathrm{3D}\text{-}\mathrm{VAR}} := \boldsymbol{\Sigma}\boldsymbol{\Sigma}^\top$.  The weight vector $\pmb{w}$ gives the linear combination of the columns of the matrix factor that describe $\pmb{x}_k$ as a perturbation.  One can iteratively optimise the 3D-VAR cost function as such with a locally quadratic approximation.  Let $\overline{\pmb{x}}^i_k$ be the $i$-th iterate for the proposal value at time $t_k$, defining $\overline{\pmb{x}}^0_k := \overline{\pmb{x}}_k^\mathrm{f}$.  The $i+1$-st iteration is defined (up to a subtlety pointed out by \textcite{menetrier2015}) as
\begin{align}
\overline{\pmb{x}}^{i+1}_k = \overline{\pmb{x}}^{i}_k + \boldsymbol{\Sigma}\overline{\pmb{w}}^{i+1},
\label{eq:iterative_state_def}
\end{align}
where $\overline{\pmb{w}}^{i+1}$ is the argmin of the following \emph{incremental cost function}
\begin{align}
\mathcal{J}_{\mathrm{3D}\text{-}\mathrm{VARI}}(\pmb{w}) := \frac{1}{2} \parallel \pmb{w}\parallel^2 + \frac{1}{2} \parallel \pmb{y}_k - \mathcal{H}_k\left(\overline{\pmb{x}}^i_k\right) - \mathbf{H}_k \boldsymbol{\Sigma} \pmb{w}\parallel^2_{\mathbf{R}_k}.\label{eq:weight_space_3dvar}
\end{align}
This approximate cost function follows from truncating the Taylor expansion
\begin{align}
\mathcal{H}_k(\pmb{x}) = \mathcal{H}_k\left(\overline{\pmb{x}}^i_k\right) + \mathbf{H}_k\boldsymbol{\Sigma}\pmb{w} + \mathcal{O}\left(\parallel\pmb{w} \parallel^2\right)
\end{align}
and the definition of the Mahalanobis distance.  Notice that $\overline{\pmb{w}}^{i+1}$ is defined uniquely as Eq.~\eqref{eq:weight_space_3dvar} is quadratic with respect to the weights, allowing one to use fast and adequate methods such as the conjugate gradient \parencite[see Chapter 5]{nocedal2006numerical}.  The iterations are set to terminate when $\parallel \overline{\pmb{w}}^i\parallel$ is sufficiently small, representing a negligible change from the last proposal. This approach to nonlinear optimisation is known as the incremental approach, and corresponds mathematically to the Gauss-Newton optimisation method, an efficient simplification of Newton's descent method for nonlinear least-squares \parencite[see Chapter 10]{nocedal2006numerical}. As described above, it consists of an iterative linearisation of the nonlinear optimisation problem, which is a concept utilised throughout this tutorial.

The 3D-VAR approach described above played an important role in early DA methodology, and aspects of this technique are still widely used in ensemble-based DA when using the regularisation technique known as covariance hybridisation \parencite{hamill2000, lorenc2003, penny2017}.  However, optimising the state with the climatological-posterior background neglects important information in the time-dependent features of the forecast spread, described sometimes as the `errors of the day' \parencite{corrazza2003}.  The \emph{extended Kalman filter} (EKF) can be interpreted to generalise the 3D-VAR cost function by including a time-dependent background, like the KF, by approximating its time-evolution at first order  with the tangent linear model as follows.

Suppose that the equations of motion are generated by a nonlinear function, independent of time for simplicity,
\begin{align}
\frac{\mathrm{d}}{\mathrm{d}t} \pmb{x} := \pmb{f}(\pmb{x}), & &\pmb{x}_k= \mathcal{M}_k(\pmb{x}_{k-1}):= \int_{t_{k-1}}^{t_k} \! \pmb{f}(\pmb{x}) \, \mathrm{d}t + \pmb{x}_{k-1}.
\end{align}
One can extend the linear-Gaussian approximation for the forecast pdf by modelling the state as a perturbation of the analysis mean,
\begin{subequations}
\begin{align}
&\pmb{x}_{k-1} := \overline{\pmb{x}}^\mathrm{a}_{k-1} + \pmb{\delta}_{k-1} \sim  \mathcal{N}\left( \overline{\pmb{x}}^\mathrm{a}_{k-1}, \mathbf{B}^\mathrm{a}_{k-1}\right),\\
\Leftrightarrow \,\,\, &\pmb{\delta}_{k-1} \sim \mathcal{N}\left(\pmb{0}, \mathbf{B}^\mathrm{a}_{k-1}\right) .
\end{align}
\end{subequations}
The evolution of the perturbation $\pmb{\delta}_{k-1}$ is written via Taylor's theorem as
\begin{subequations}
\begin{align}
\frac{\mathrm{d}}{\mathrm{d}t} \pmb{\delta}_{k-1}&:= \frac{\mathrm{d}}{\mathrm{d}t}\left( \pmb{x}_{k-1} - \overline{\pmb{x}}_{k-1}^\mathrm{a}\right)\\
&=\pmb{f}\left(\pmb{x}_{k-1}\right) - \pmb{f}\left(\overline{\pmb{x}}_{k-1}^\mathrm{a}\right)\\
&=\nabla_{\pmb{x}}\pmb{f}\left(\overline{\pmb{x}}_{k-1}^\mathrm{a}\right)\pmb{\delta}_{k-1} + \mathcal{O}\left(\parallel \pmb{\delta}_{k-1}\parallel^2\right),
\end{align}
\end{subequations}
where $\nabla_{\pmb{x}}\pmb{f}(\overline{\pmb{x}}_{k-1}^\mathrm{a})$ is the Jacobian equation with dependence on the underlying analysis mean state.  The linear evolution defined by the truncated Taylor expansion about the underlying reference solution $\overline{\pmb{x}}$,
\begin{align}
\frac{\mathrm{d}}{\mathrm{d}t} \pmb{\delta}:= \nabla_{\pmb{x}}\pmb{f}(\overline{\pmb{x}})\cdot \pmb{\delta},
\label{eq:tangent_linear_model}
\end{align}
is known as the \emph{tangent linear model}.

Making the approximation of the tangent linear model for the first order evolution of the modelled state,
\begin{subequations}
\begin{align}
&\frac{\mathrm{d}}{\mathrm{d}t} \pmb{x} \approx \pmb{f}(\overline{\pmb{x}}) + \nabla_{\pmb{x}}\pmb{f}(\overline{\pmb{x}})\cdot \pmb{\delta}\\
\Rightarrow&\int_{t_{k-1}}^{t_k} \! \frac{\mathrm{d}}{\mathrm{d}t}\pmb{x} \, \mathrm{d}t \approx  \int_{t_{k-1}}^{t_k} \! \pmb{f}(\overline{\pmb{x}}) \, \mathrm{d}t + \int_{t_{k-1}}^{t_k} \! \nabla_{\pmb{x}}\pmb{f}(\overline{\pmb{x}})\cdot \pmb{\delta} \, \mathrm{d}t \\
\Rightarrow & \pmb{x}_{k} \approx \mathcal{M}_{k}\left(\overline{\pmb{x}}_{k-1}\right) + \mathbf{M}_k\pmb{\delta}_{k-1}
\end{align}
\end{subequations}
where $\mathbf{M}_k$ is the resolvent of the tangent linear model.  Given that Gaussians are closed under affine transformations, the EKF approximation for the (perfect) evolution of the state vector is defined as
\begin{align}
p(\pmb{x}_{k}|\pmb{y}_{k-1:1}) \approx n\left(\pmb{x}_k |\mathcal{M}_k\left(\overline{\pmb{x}}^\mathrm{a}_{k-1}\right), \mathbf{M}_k \mathbf{B}^\mathrm{a}_{k-1}\mathbf{M}^\top_{k}\right).
\end{align}
Respectively, define the EKF analysis cost function as
\begin{align}
\mathcal{J}_{\mathrm{EKF}}(\pmb{x}) := \frac{1}{2}\parallel \overline{\pmb{x}}_k^\mathrm{f}  - \pmb{x} \parallel_{\mathbf{B}^\mathrm{f}_{k}}^2 + \frac{1}{2} \parallel \pmb{y}_k - \mathcal{H}_k\left(\pmb{x}\right)\parallel^2_{\mathbf{R}_k},
\end{align}
where
\begin{subequations}
\begin{align}
\overline{\pmb{x}}_k^\mathrm{f}&:=  \mathcal{M}_k\left(\overline{\pmb{x}}^\mathrm{a}_{k-1}\right),\\
\mathbf{B}_k^\mathrm{f} &:=\mathbf{M}_k \mathbf{B}^\mathrm{a}_{k-1}\mathbf{M}^\top_{k}.
\end{align}
\end{subequations}
A similar, locally quadratic, weight-space optimisation can be performed versus the nonlinear observation operator using a matrix decomposition for the time-varying background forecast error covariance, defining the state as a perturbation using this matrix factor similar to Eq.~\eqref{eq:weight_definition}.  This gives the square root EKF, which was used historically to improve the stability over the direct approach \parencite{tippett2003}.

The accuracy and stability of the above EKF depends strongly on the length of the forecast horizon $\Delta_t$ \parencite{miller1994}. For short-range forecasting, the perturbation dynamics of the tangent linear model can be an adequate approximation, and is an underlying approximation for most operational DA \parencite{carrassi2018}.  However, the explicit linear approximation can degrade quickly, and especially when the mean state is not accurately known.  The perturbation dynamics in the linearisation about the approximate mean can differ substantially from the nonlinear dynamics of the true system as the approximate mean state diverges from the modelled system \parencite[see Chapters 7 - 8]{grewal2014kalman}.  These stability issues, and the computational cost of explicitly representing the evolution of the background covariance in the tangent linear model, limits the use of the EKF for geophysical DA.  Nonetheless, this technique provides important intuition that is developed later in this tutorial.

\subsection{The Gaussian 4D Cost Function}
Consider again the perfect, linear-Gaussian model  represented by Eq.~\eqref{eq:linear_gaussian}.  Using the Markov assumption, Eq.~\eqref{eq:markov_assumption}, and the independence of observation errors, Eq.~\eqref{eq:observation_independence}, recursively for the hidden Markov model, the joint posterior decomposes as 
\begin{align}
p(\pmb{x}_{L:0} \vert \pmb{y}_{L:1}) &\propto  \underbrace{ p(\pmb{x}_0) }_{(i)}
\underbrace{ \left[\prod_{k=1}^L p(\pmb{x}_k|\pmb{x}_{k-1}) \right]}_{(ii)}
\underbrace{ \left[ \prod_{k=1}^L   p(\pmb{y}_k \vert \pmb{x}_k) \right]}_{(iii)} \label{eq:4d_bayes}
\end{align}
where (i) is the prior for the initial state $\pmb{x}_0$; (ii) is the free-forecast with the perfect model $\mathbf{M}_k$ depending on some outcome drawn from the prior; and (iii) is the joint likelihood of the observations in the DAW, given the background forecast.

Define the composition of the linear model forecast from time $t_{k-1}$ to $t_l$ as
\begin{align}
\mathbf{M}_{l:k} := \mathbf{M}_{l}\cdots \mathbf{M}_{k}, & & \mathbf{M}_{k:k} := \mathbf{I}_{N_x}.
\end{align}
Using the perfect, linear model hypothesis, note that
\begin{align}
\pmb{x}_k := \mathbf{M}_{k:1}\pmb{x}_{0} 
\end{align}
for every $k$.  Therefore, the transition pdfs in Eq.~\eqref{eq:4d_bayes} are reduced to a trivial condition by re-writing
\begin{subequations}
\begin{align}
p(\pmb{x}_{L:0} \vert \pmb{y}_{L:1}) &\propto  p(\pmb{x}_0)  \left[\prod_{k=1}^L p(\pmb{x}_k|\mathbf{M}_{k-1:1}\pmb{x}_{0}) \right]\left[ \prod_{k=1}^L   p(\pmb{y}_k \vert \pmb{x}_k) \right] \\
&\propto p(\pmb{x}_0) \left[ \prod_{k=1}^L   p(\pmb{y}_k \vert \mathbf{M}_{k:1}\pmb{x}_0) \right]
\end{align}
\end{subequations}
as this pdf evaluates to zero whenever $\pmb{x}_k \neq \mathbf{M}_{k:1}\pmb{x}_0$.  Given a Gaussian prior, and Gaussian observation error distributions, the minus-log-posterior \emph{four-dimensional} (4D) cost function is derived as
\begin{align}
\mathcal{J}_{\mathrm{4D}}(\pmb{x}_0) :=  \frac{1}{2}\parallel \overline{\pmb{x}}_0 - \pmb{x}_0\parallel^2_{\mathbf{B}_0} + \frac{1}{2}\sum_{k=1}^L \parallel \pmb{y}_k - \mathbf{H}_k \mathbf{M}_{k:1}\pmb{x}_0\parallel_{\mathbf{R}_k}^2.
\label{eq:4d_map_cost}
\end{align}

Notice that the 4D cost function in Eq.~\eqref{eq:4d_map_cost} is actually quadratic with respect to the initial condition $\pmb{x}_0$. Therefore, the smoothing problem in the perfect, linear-Gaussian model has a unique, optimal initial condition that minimises the sum-of-square deviations from the prior mean, with distance weighted with respect to the prior covariance, and the observations, with distance weighted with respect to the observation error covariances.  The optimal, smoothed initial condition using observation information up to time $t_L$ is denoted $\overline{\pmb{x}}^\mathrm{s}_{0|L}$; this gives the smoothed model states at subsequent times by the perfect model evolution,
\begin{align}
\overline{\pmb{x}}^\mathrm{s}_{k|L} := \mathbf{M}_{k:1} \overline{\pmb{x}}^\mathrm{s}_{0|L}.
\end{align}
The derivation above is formulated for a smoothing problem with a static-in-time DAW, in which an entire time series of observations $\pmb{y}_{L:1}$ is available, and for which the estimation may be performed `offline'.  However, a simple extension of this analysis accommodates DAWs that are sequentially shifted in time, allowing an `online' smoothing estimate to be formed analogously to sequential filtering.  

\begin{figure*}
\center
\includegraphics[width=\linewidth]{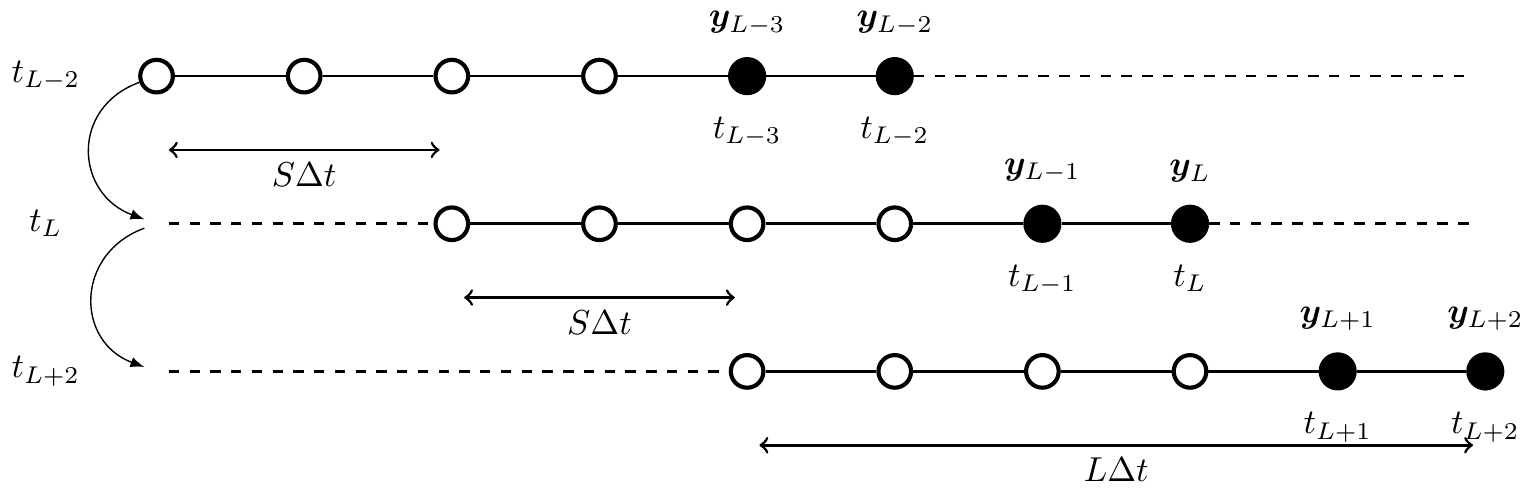}
\caption{Three cycles of a shift $S=2$, lag $L=5$ fixed-lag smoother, cycle number is increasing top to bottom.  Time indices on the left-hand margin indicate the current time for the associated cycle of the algorithm.  New observations entering the current DAW are shaded black.  Source: \textcite{grudzien2021fast}, adapted from \textcite{asch2016}.}
\label{fig:lag_shift}
\end{figure*}

Fixed-lag smoothing sets the length of the DAW, $L$, to be fixed for all time, while the underlying states are cycled through this window as it shifts forwards in time.  Given the \emph{lag} of length $L$, suppose that a \emph{shift} $S$ is defined for which $1 \leq S\leq L$.  It is convenient to consider an algorithmically stationary DAW, referring to the time indices $\{t_1, \cdots, t_L\}$.   In a given cycle, the joint posterior $p(\pmb{x}_{L:1}|\pmb{y}_{L:1})$ is estimated.  After the estimate is produced, the DAW is subsequently shifted in time by $S \times \Delta_t$, and all states are re-indexed by $t_k \longleftarrow t_{k+S}$ to begin the next cycle.  For a lag of $L$ and a shift of $S$, the observation vectors at times $\{t_{L-S+1},\cdots, t_{L}\}$ correspond to observations newly entering the DAW for the analysis performed up to time $t_L$.  When $S=L$, the DAWs are disconnected and adjacent in time, whereas for $S<L$ there is an overlap between the estimated states in sequential DAWs.  Figure \ref{fig:lag_shift} provides a schematic of how the DAW is shifted for a lag of $L=5$ and shift $S=2$.   Following the common convention in DA that there is no observation at time zero, in addition to the DAW, $\{t_1,\cdots, t_L\}$, states at time $t_0$ may be estimated or utilised in order to connect estimates between adjacent / overlapping DAWs. 

Consider the algorithmically stationary DAW, and suppose that the current analysis time is $t_L$, where the joint posterior pdf $p(\pmb{x}_{L-S:1-S} | \pmb{y}_{L-S: 1-S})$ is available from the last fixed-lag smoothing cycle at analysis time $t_{L-S}$.  Using the independence of observation errors and the Markov assumption recursively,
\begin{align}
\begin{split}
p(\pmb{x}_{L:1} & \vert \pmb{y}_{L:1-S}) \propto \\
&\int \! \underbrace{\mathrm{d}\pmb{x}_{0} \, p\left(\pmb{x}_0 | \pmb{y}_{L-S:1-S}\right) }_{(i)} \underbrace{\left[\prod_{k=1}^{L} p(\pmb{x}_k \vert \pmb{x}_{k-1})\right]}_{(ii)} \underbrace{\left[ \prod_{k=L-S+1}^L p(\pmb{y}_k \vert \pmb{x}_k )\right]}_{(iii)}
,
\end{split}
\end{align}
where: (i) now represents averaging out the initial condition at time $t_0$ with respect to the marginal smoothing pdf for $\pmb{x}^\mathrm{s}_{0|L-S}$ over the last DAW; (ii) represents the free forecast of the smoothed estimate for $\pmb{x}^\mathrm{s}_{0|L-S}$; and term (iii) represents the joint likelihood of the newly incoming observations to the DAW given the forecasted model state.  Noting that $p(\pmb{x}_{L:1}| \pmb{y}_{L:1}) \propto p(\pmb{x}_{L:1} | \pmb{y}_{L:1-S})$, this provides a recursive form of the 4D cost function, shifted sequentially in time,
\begin{align}
\mathcal{J}_{\mathrm{4D}\text{-}\mathrm{seq}}(\pmb{x}_0) :=  \frac{1}{2}\parallel \overline{\pmb{x}}_{0|L-S}^\mathrm{s} - \pmb{x}_0\parallel^2_{\mathbf{B}_{0|L-S}^\mathrm{s}} + \frac{1}{2}\sum_{k=L-S+1}^L \parallel \pmb{y}_k - \mathbf{H}_k \mathbf{M}_{k:1}\pmb{x}_0\parallel_{\mathbf{R}_k}^2,
\label{eq:4d_map_cost_sequential}
\end{align}
where $\mathbf{B}_{0|L-S}^\mathrm{s}$ refers to the smoothed error covariance from the last DA cycle.  As with the KF cost function, the posterior error covariance, conditioning on observations up to time $t_L$, is identified with
\begin{align}
\mathbf{B}^\mathrm{s}_{0|L}:= \boldsymbol{\Xi}_{\mathcal{J}_{\mathrm{4D}\text{-}\mathrm{seq}}}^{-1}.
\end{align}
Given the perfect, linear-Gaussian model assumption, the mean and covariance are propagated to time $t_S$ via
\begin{subequations}
\begin{align}
\overline{\pmb{x}}^\mathrm{s}_{S|L} &:= \mathbf{M}_{S:1}\overline{\pmb{x}}_{0|L}^\mathrm{s},\\
\mathbf{B}_{S|L}^\mathrm{s} &:= \mathbf{M}_{S:1} \mathbf{B}_{0|L}^\mathrm{s} \mathbf{M}_{S:1}^\top,
\end{align}
\end{subequations}
and states are re-indexed as $t_{k+S}\longleftarrow t_k$ to initialise the next cycle.  This provides a 4D derivation of the \emph{Kalman smoother} (KS) assuming a perfect model, though not all developments take the above approach, using a global analysis over all observations in the current DAW at once.  Other formulations use an alternating: (i) forwards sequential filtering pass to update the current forecast; and (ii) a backwards-in-time sequential filtering pass over the DAW to condition lagged states on the new information. Examples of smoothers that follow this alternating analysis include the \emph{ensemble Kalman smoother} (EnKS) \parencite{evensen2000} and the ensemble \emph{Rauch-Tung-Striebel} (RTS) smoother \parencite{raanes2016}, with the EnKS to be discussed in Section \ref{section:enks}.

\subsection{Incremental 4D-VAR}
The method of incremental \emph{four-dimensional, variational} (4D-VAR) data assimilation \parencite{ledimet1986,talagrand1987,courtier1994} is a classical and widely used DA technique that extends the linear-Gaussian 4D analysis to nonlinear settings, both for fixed and sequential DAWs.  Modern formulations of the 4D-VAR analysis furthermore handle model errors, as in weak-constraint 4D-VAR \parencite{tremolet2006,desroziers2014,laloyaux2020}, and may include time-varying background error covariances by combining an ensemble of 4D-VAR \parencite{bonavita2012}, which is known as \emph{ensemble of data assimilation} (EDA).

Assuming now that the state and observation models are nonlinear, as in Eq.~\eqref{eq:perfect_model}, denote the composition of the nonlinear model forecast from time $t_{k-1}$ to time $t_{l}$ as
\begin{align}
\mathcal{M}_{l:k}:= \mathcal{M}_{l} \circ \cdots \circ \mathcal{M}_{k}, & & \mathcal{M}_{k:k}:= \mathbf{I}_{N_x}.
\end{align}
Then, the linear-Gaussian 4D cost function is formally extended as a Gaussian approximation in 4D-VAR with
\begin{align}
\mathcal{J}_{\mathrm{4D}\text{-}\mathrm{VAR}}(\pmb{x}_0) :=  \frac{1}{2}\parallel \overline{\pmb{x}}_0 - \pmb{x}_0\parallel^2_{\mathbf{B}_{\mathrm{4D}\text{-}\mathrm{VAR}}} + \frac{1}{2}\sum_{k=1}^L \parallel \pmb{y}_k - \mathcal{H}_k\circ \mathcal{M}_{k:1}(\pmb{x}_0)\parallel_{\mathbf{R}_k}^2,
\label{eq:4d_var_primal}
\end{align}
where $\mathbf{B}_{\mathrm{4D}\text{-}\mathrm{VAR}} \longleftarrow \mathbf{B}_{\mathrm{3D}\text{-}\mathrm{VAR}}$.  A similar indexing to equation to Eq.~\eqref{eq:4d_map_cost_sequential} gives the sequential form over newly incoming observations.

The incremental linearisation that was performed in 3D-VAR in Eq.~\eqref{eq:weight_space_3dvar} forms the basis for the classical technique of incremental 4D-VAR.  Suppose again that an explicit Cholesky factor for the background covariance is given $\mathbf{B}_{\mathrm{4D}\text{-}\mathrm{VAR}}:= \boldsymbol{\Sigma}\boldsymbol{\Sigma}^\top$, where the state is written as a perturbation of the $i$-th iterate as
\begin{align}
\pmb{x}_0: = \overline{\pmb{x}}^i_0 + \boldsymbol{\Sigma} \pmb{w}.
\end{align}
Taking a Taylor expansion of the composition of the nonlinear observation model and the nonlinear state model at the $i$-th iterate,
\begin{align}
\mathcal{H}_k \circ \mathcal{M}_{k:1}\left(\pmb{x}_0\right) = \mathcal{H}\circ\mathcal{M}_{k:1}\left(\overline{\pmb{x}}^i_0\right) + \mathbf{H}_k \mathbf{M}_{k:1} \boldsymbol{\Sigma} \pmb{w} + \mathcal{O}\left(\parallel \pmb{w} \parallel^2\right),
\end{align}
the incremental cost function is rendered
\begin{align}
\mathcal{J}_{\mathrm{4D}\text{-}\mathrm{VARI}}(\pmb{w}) = \frac{1}{2} \parallel \pmb{w}\parallel^2 + \sum_{k=1}^L \parallel \pmb{y}_k - \mathcal{H}_k \circ \mathcal{M}_{k:1}\left( \overline{\pmb{x}}^i_0 \right) - \mathbf{H}_k\mathbf{M}_{k:1} \boldsymbol{\Sigma} \pmb{w}\parallel_{\mathbf{R}_k}^2.
\label{eq:incremental_4dvar}
\end{align}
The argmin of Eq.~\eqref{eq:incremental_4dvar} is defined as $\overline{\pmb{w}}^i$, where $\overline{\pmb{x}}^{i+1}$ is defined as in Eq.~\eqref{eq:iterative_state_def}.

It is important to remember that in Eq.~\eqref{eq:incremental_4dvar}, the terms $\mathbf{M}_{k:1}$ involve the calculation of the tangent linear model with reference to the underlying nonlinear solution $\mathcal{M}_{L:1}\left(\overline{\pmb{x}}^i_0\right)$, simulated freely over the entire DAW.  A key aspect to the efficiency of the incremental 4D-VAR approach is in the adjoint model computation of the gradient to this cost function, commonly known as \emph{backpropagation} in statistical and machine learning \parencite[see Chapter 7]{rojas2013neural}.    The \emph{adjoint model} with respect to the proposal is defined as
\begin{align}
\frac{\mathrm{d}}{\mathrm{d}t} \tilde{\pmb{\delta}} = \left(\nabla_{\pmb{x}} \pmb{f}(\overline{\pmb{x}}^i)\right)^\top \tilde{\pmb{\delta}}.
\label{eq:adjoint_model}
\end{align}
The solution to the linear adjoint model above defines the adjoint resolvent matrix $\mathbf{M}_k^\top$, i.e., the transpose of the tangent linear resolvent.  Notice that for the composition of the tangent linear model forecasts with $k\leq l$, the adjoint is given as
\begin{align}
\mathbf{M}_{l:k}^\top := \left(\mathbf{M}_l \mathbf{M}_{l-1} \cdots \mathbf{M}_{k+1} \mathbf{M}_k \right)^\top = \mathbf{M}_k^\top \mathbf{M}_{k+1}^\top \cdots \mathbf{M}_{l-1}^\top \mathbf{M}_l^\top.
\end{align}
This means that the adjoint model variables $\tilde{\pmb{\delta}}_l$ are propagated by the linear resolvents of the adjoint model, but applied in reverse chronological order from the tangent linear model from last time $t_l$ to the initial time $t_{k-1}$, 
\begin{align}
\tilde{\pmb{\delta}}_{k-1} := \mathbf{M}^\top_{l:k} \tilde{\pmb{\delta}}_l
\end{align}
transmitting the sensitivity from a future time back to a perturbation of $\overline{\pmb{x}}_{k-1}^i$ \parencite[see Section 6.3]{kalnay2003}.

Define the \emph{innovation} vector of the $i$-th iterate and the $k$-th observation vector as
\begin{align}
\overline{\pmb{\delta}}^i_k := \pmb{y}_k - \mathcal{H}_k \circ \mathcal{M}_{k:1}\left( \overline{\pmb{x}}^i_0 \right) .
\label{eq:innovation}
\end{align}
Then the gradient of Eq.~\eqref{eq:incremental_4dvar} with respect to the weight vector is written
\begin{align}
\nabla_{\pmb{w}} \mathcal{J}_{\mathrm{4D}\text{-}\mathrm{VARI}} = \pmb{w} - \sum_{k=1}^L \boldsymbol{\Sigma}^\top \mathbf{M}_{k:1}^\top \mathbf{H}_k^\top\mathbf{R}_k^{-1} \left[ \overline{\pmb{\delta}}^i_{k} - \mathbf{H}_k \mathbf{M}_{k:1}\boldsymbol{\Sigma} \pmb{w} \right].
\end{align}
Making the substitution
\begin{align}
\pmb{\Delta}_k:= \mathbf{R}_k^{-1} \left[ \overline{\pmb{\delta}}^i_{k} - \mathbf{H}_k \mathbf{M}_{k:1}\boldsymbol{\Sigma} \pmb{w} \right],
\label{eq:obs_pert}
\end{align}
notice that the gradient is written in Horner factorisation as
\begin{align}
\nabla_{\pmb{w}} \mathcal{J}_{\mathrm{4D}\text{-}\mathrm{VARI}} = \pmb{w} - \boldsymbol{\Sigma}^\top \mathbf{M}_1^\top \left[ \mathbf{H}^\top_1 \pmb{\Delta}_1 + \mathbf{M}^\top_2 \left[ \mathbf{H}^\top_2 \pmb{\Delta}_2 +\cdots + \left[\mathbf{M}_L^\top\mathbf{H}^\top_L \pmb{\Delta}_L \right] \right]\right] .
\end{align}
Defining the adjoint variable as $\tilde{\pmb{\delta}}_L:= \mathbf{H}_L^\top \pmb{\Delta}_L$, with recursion in $t_k$,
\begin{align}
\tilde{\pmb{\delta}}_k &:= \mathbf{H}_k^\top \pmb{\Delta}_k + \mathbf{M}_{k+1}^\top \tilde{\pmb{\delta}}_{k+1},
\label{eq:discrete_adjoint}
\end{align}
the gradient of the cost function is written as
\begin{align}
\nabla_{\pmb{w}} \mathcal{J}_{\mathrm{4D}\text{-}\mathrm{VARI}}:=\pmb{w} - \boldsymbol{\Sigma}\tilde{\pmb{\delta}}_0.
\label{eq:adjoint_gradient}
\end{align}
Using the above definitions, the gradient of the incremental cost function is computed from the following steps: (i) a forwards pass of the free evolution of $\pmb{x}_0^i$ under the nonlinear forecast model, computing the innovations as in Eq.~\eqref{eq:innovation}; (ii) the propagation of the perturbation $\boldsymbol{\Sigma}\pmb{w}$ in the tangent linear model and the linearised observation operator with respect to the proposal, $\overline{\pmb{x}}^i_k$, in order to compute the terms in Eq.~\eqref{eq:obs_pert}; and (iii) the back propagation of the sensitivities in the adjoint model, Eq.~\eqref{eq:adjoint_model},
to obtain the adjoint variables recursively as in Eq.~\eqref{eq:discrete_adjoint} and the gradient from Eq.~\eqref{eq:adjoint_gradient}.  

The benefit of the above approach is that this gives an extremely efficient calculation of the gradient, provided that the tangent linear and adjoint models are available, which, in turn, is key to very efficient numerical optimisations of cost functions. However, the trade off is that the tangent linear and adjoint models of a full scale geophysical state model, and of the observation model, require considerable development time and expertise.  Increasingly, these models can be computed abstractly by differentiating a computer program alone, in what is known as \emph{automatic differentiation} of code \parencite{griewank1989, griewank2003,hascoet2014,baydin2018}.  When tangent linear and adjoint models are not available, one alternative is to use ensemble sampling techniques in the fully nonlinear model $\mathcal{M}$ alone.  The ensemble-based analysis provides a complementary approach to the explicit tangent linear and adjoint model analysis -- this approach can be developed independently, or hybridised with the use of the tangent linear and adjoint models as in various flavours of hybrid \emph{ensemble-variational} (EnVAR) techniques \parencite{asch2016, bannister2017}.  This tutorial focuses on the independent development of EnVAR estimators in the following section.

\section{Bayesian Ensemble-Variational Estimators}

\subsection{The Ensemble Transform Kalman Filter}
Consider once again the perfect, linear-Gaussian model in Eq.~\eqref{eq:linear_gaussian}.  Rather than explicitly computing the evolution of the background mean and error covariance in the linear model as in Eq.~\eqref{eq:kf_forecast} and the KF, one can alternatively estimate the state mean and the error covariances using a statistical sampling approach.  Let $\left\{\pmb{x}_{k,i}^{\mathrm{f/a}}\right\}_{i=1}^{N_e}$ be replicates of the model state, \emph{independently and identically distributed} (iid) according to the distribution
\begin{align}
\pmb{x}_{k,i}^{\mathrm{f/a}} \sim \mathcal{N}\left(\overline{\pmb{x}}^\mathrm{f/a}_k, \mathbf{B}_k^\mathrm{f/a}\right).
\label{eq:columns_distributed}
\end{align}
Given the iid assumption above, the ensemble-based mean and the ensemble-based covariance
\begin{subequations}
\begin{align}
\hat{\pmb{x}}^\mathrm{f/a}_k &:= \frac{1}{N_e}\sum_{i=1}^{N_e} \pmb{x}_{k,i}^{\mathrm{f/a}}, \\
\mathbf{P}_k^\mathrm{f/a} &:= \frac{1}{N_e -1} \sum_{i=1}^{N_e} \left(\pmb{x}_{k,i}^{\mathrm{f/a}} - \hat{\pmb{x}}^\mathrm{f/a}_k\right)\left(\pmb{x}_{k,i}^{\mathrm{f/a}} - \hat{\pmb{x}}^\mathrm{f/a}_k\right)^\top,
\end{align}
\label{eq:series_sample_stats}
\end{subequations}
are (asymptotically) consistent estimators of the background, i.e.,
\begin{align}
\mathbb{E}\left[ \hat{\pmb{x}}^\mathrm{f/a}_k\right]= \overline{\pmb{x}}^\mathrm{f/a}_k, & &
\lim_{N_e \rightarrow \infty}\mathbb{E}\left[\mathbf{P}_k^\mathrm{f/a} \right]= \mathbf{B}_k^{\mathrm{f/a}},
\end{align}
where the expectation is over the possible realisations of the random sample.
Particularly, the multivariate central limit theorem gives that
\begin{align}
\left(\mathbf{P}_k^\mathrm{f/a}\right)^{-\frac{1}{2}} \left( \pmb{x}^\mathrm{f/a}_k - \hat{\pmb{x}}^\mathrm{f/a}_k \right) \rightarrow_d \mathcal{N}\left(\pmb{0}, \mathbf{I}_{N_x}\right),\label{eq:multivariate_clt}
\end{align}
referring to convergence in distribution as $N_e\rightarrow \infty$ \parencite[see Section 6.2]{hardle2017basic}.  

Using the relationships in Eq.~\eqref{eq:series_sample_stats}, these estimators are efficiently encoded as linear operations on the ensemble matrix.  Define the \emph{ensemble matrix} and the \emph{perturbation matrix} as
\begin{align}
\mathbf{E}^\mathrm{f/a}_k&:= \begin{pmatrix} \pmb{x}_{k,1}^{\mathrm{f/a}} &\cdots & \pmb{x}_{k,N_e}^{\mathrm{f/a}}\end{pmatrix} \in \mathbb{R}^{N_x \times N_e},\\
\mathbf{X}^\mathrm{f/a}_k&:= \begin{pmatrix} \pmb{x}_{k,1}^{\mathrm{f/a}} -\hat{\pmb{x}}_k^\mathrm{f/a}&\cdots & \pmb{x}_{k,N_e}^{\mathrm{f/a}} -\hat{\pmb{x}}_k^\mathrm{f/a}\end{pmatrix} \in \mathbb{R}^{N_x \times N_e},
\end{align}
i.e., as the arrays with the columns given as the ordered replicates of the model state and their deviations from the ensemble mean respectively.  For $\pmb{1}$ defined as the vector composed entirely of ones, define the following linear operations conformally in their dimensions
\begin{subequations}
\begin{align}
\hat{\pmb{x}}_k^\mathrm{f/a} &= \mathbf{E}^\mathrm{f/a}_k \pmb{1}/N_e, \\
\mathbf{X}_k^\mathrm{f/a} &= \mathbf{E}_k^\mathrm{f/a} \left(\mathbf{I}_{N_e} - \pmb{1}\pmb{1}^\top/N_e \right),\\
\mathbf{P}_k^\mathrm{f/a} &= \left(\mathbf{X}_k^\mathrm{f/a} \right)\left(\mathbf{X}_k^\mathrm{f/a} \right)^\top / (N_e -1).
\end{align}
\end{subequations}
The operator $\mathbf{I}_{N_e} - \pmb{1}\pmb{1}^\top/N_e$ is the orthogonal complementary projection operator to the span of the vector of ones, known as the \emph{centring operator} in statistics \parencite[see Section 6.1]{hardle2017basic}; this has the effect of transforming the ensemble to mean zero.

Notice, when $N_e \leq N_x$ the ensemble-based error covariance has rank of at most $N_e -1 $ irrespective of the rank of the background error covariance, corresponding to the one degree of freedom lost in computing the ensemble mean.  Therefore, to utilise the ensemble error covariance in the least-squares optimisation as with the KF in Eq.~\eqref{eq:simple_map}, a new construction is necessary.  For a generic matrix $\mathbf{A}\in \mathbb{R}^{N\times M}$ with full column rank $M$, define the (left Moore-Penrose) \emph{pseudo-inverse} \parencite[see page 423]{meyer2000}
\begin{align}
\mathbf{A}^\dagger := \left(\mathbf{A}^\top \mathbf{A}\right)^{-1}\mathbf{A}^\top .
\end{align}
In particular, $\mathbf{A}^\dagger \mathbf{A} = \mathbf{I}_M$ and the orthogonal projector into the column span of $\mathbf{A}$ is defined by $\mathbf{A}\mathbf{A}^\dagger$.  When $\mathbf{A}$ has full column rank as above, define the Mahalanobis `distance' with respect to $\mathbf{G} = \mathbf{A}\mathbf{A}^\top$ as
\begin{align}
\parallel \pmb{v} \parallel_{\mathbf{G}} := \sqrt{ \left(\mathbf{A}^\dagger\pmb{v}\right)^\top \left(\mathbf{A}^\dagger \pmb{v}\right)}.
\end{align}
Note that in the case that $\mathbf{G}$ does not have full column rank, i.e., $N>M$, this is not a true norm on $\mathbb{R}^N$ as it is degenerate in the null space of $\mathbf{A}^\dagger$.  This instead represents a lift of a non-degenerate norm in the column span of $\mathbf{A}$ to $\mathbf{R}^N$.  In the case that $\pmb{v}$ is in the column span of $\mathbf{A}$, one can equivalently write
\begin{subequations}
\begin{align}
\pmb{v} = \mathbf{A}\pmb{w},  \\
\parallel \pmb{v}\parallel_\mathbf{G} = \parallel \pmb{w}\parallel,
\end{align}
\end{subequations}
for a vector of weights $\pmb{w}\in \mathbf{R}^M$.

The \emph{ensemble Kalman filter} (EnKF) cost function for the linear-Gaussian model is defined
\begin{subequations}
\begin{align}
& \mathcal{J}_{\mathrm{EnKF}}(\pmb{x}_k):= \frac{1}{2}\parallel \hat{\pmb{x}}_k^\mathrm{f} - \pmb{x}_k \parallel_{\mathbf{P}_k^\mathrm{f}}^2 + \frac{1}{2}\parallel \pmb{y}_k - \mathbf{H}_k \pmb{x}_k \parallel_{\mathbf{R}_k}^2
\label{eq:ensemble_cost_state_space}
\\
\Leftrightarrow \,\, & \mathcal{J}_{\mathrm{EnKF}}(\pmb{w}) := \frac{1}{2}(N_e - 1)\parallel \pmb{w}\parallel^2 +\frac{1}{2}\parallel \pmb{y}_k - \mathbf{H}_k\hat{\pmb{x}}^\mathrm{f}_k  - \mathbf{H}_k \mathbf{X}_k^\mathrm{f} \pmb{w}\parallel_{\mathbf{R}_k}^2 ,
\label{eq:ensemble_cost_weight_space}
\end{align}
\end{subequations}
where the model state is written as a perturbation of the ensemble mean
\begin{align}
\pmb{x}_k = \hat{\pmb{x}}_k^\mathrm{f} + \mathbf{X}_k^\mathrm{f} \pmb{w}.
\end{align}
Notice, $\pmb{w}\in \mathbb{R}^{N_e}$, giving the linear combination of the ensemble perturbations used to represent the model state.

Define $\hat{\pmb{w}}$ to be the argmin of the cost function in Eq.~\eqref{eq:ensemble_cost_weight_space}.  \textcite{hunt2007} and \textcite{bocquet2011} demonstrate that, up to a gauge transformation, $\hat{\pmb{w}}$ yields the argmin of the state-space cost function, Eq.~\eqref{eq:ensemble_cost_state_space}, when the estimate is restricted to the ensemble span.  Equation \eqref{eq:ensemble_cost_weight_space} is quadratic in $\pmb{w}$ and can be solved to render
\begin{subequations}
\begin{align}
\hat{\pmb{w}} &:= \pmb{0} - \boldsymbol{\Xi}^{-1}_{\mathcal{J}_\mathrm{EnKF}} \nabla_{\pmb{w}} \mathcal{J}_{\mathrm{EnKF}}|_{\pmb{w}=\pmb{0}},\label{eq:enkf_newton}\\
\mathbf{T} &:= \boldsymbol{\Xi}^{-\frac{1}{2}}_{\mathcal{J}_{\mathrm{EnKF}}} ,\label{eq:right_transform}\\
\mathbf{P}_k^\mathrm{a} &= \left(\mathbf{X}_k^\mathrm{f}\mathbf{T}\right)  \left(\mathbf{X}_k^\mathrm{f} \mathbf{T}\right)^\top/(N_e -1),
\end{align}
\end{subequations}
corresponding to a single iteration of Newton's descent algorithm in Eq.~\eqref{eq:enkf_newton}, initialised with the ensemble mean, to find the optimal weights. 

The linear \emph{ensemble transform Kalman filter} (ETKF) equations are then given by
\parencite{bishop2001,hunt2007}
\begin{subequations}
\begin{align}
\mathbf{E}_k^\mathrm{f} &= \mathbf{M}_k \mathbf{E}_{k-1}^\mathrm{a} ,\\
\mathbf{E}^\mathrm{a}_k &= \hat{\pmb{x}}_k^\mathrm{f} \pmb{1}^\top + \mathbf{X}_k^\mathrm{f}\left(\hat{\pmb{w}}\pmb{1}^\top + \sqrt{N_e - 1}\mathbf{T} \mathbf{U}\right),\label{eq:ensemble_transform_part_I}
\end{align}
\end{subequations}
where $\mathbf{U}\in\mathbb{R}^{N_e \times N_e}$ can be any mean-preserving, orthogonal transformation, i.e., $\mathbf{U}\pmb{1}=\pmb{1}$.  The simple choice of $\mathbf{U} := \mathbf{I}_{N_e}$ is sufficient, but it has been demonstrated that choosing a random, mean-preserving orthogonal transformation at each analysis as above can improve the accuracy and robustness of the ETKF, smoothing out higher-order artefacts in the empirical covariance estimate \parencite{sakov2008b}. 

Notice that Eq.~\eqref{eq:ensemble_transform_part_I} is written equivalently as a single right ensemble transformation:
\begin{subequations}
\label{eq:right_ensemble_update}
\begin{align}
\mathbf{E}_k^\mathrm{a} &= \mathbf{E}_k^\mathrm{f}\boldsymbol{\Psi}_k, \\
\boldsymbol{\Psi}_k &:= \pmb{1}\pmb{1}^\top/N_e + \left(\mathbf{I}_{N_e} -\pmb{1}\pmb{1}^\top/N_e \right) \left(\hat{\pmb{w}} \pmb{1}^\top + \sqrt{N_e - 1} \mathbf{T}\mathbf{U}\right) ,
\end{align}
\end{subequations}
where the columns are approximately distributed as in Eq.~\eqref{eq:columns_distributed}, with the (asymptotic) consistency as with the central limit theorem, Eq.~\eqref{eq:multivariate_clt}.  However, in the small feasible sample sizes for realistic geophysical models, this approximation can lead to a systematic underestimation of the uncertainty of the analysis state, where the ensemble-based error covariance can become overly confident in its own estimate with artificially small variances.  Covariance inflation is a technique that is widely used to regularise the ensemble-based error covariance by increasing the empirical variances of the estimate. This can be used to handle the inaccuracy of the estimator due to the finite sample size approximation of the background mean and error covariance as in Eq.~\eqref{eq:multivariate_clt}, as well as inaccuracies due to a variety of other sources of error \parencite{carrassi2018, raanes2019adaptive, tandeo2020}.

A Bayesian hierarchical approach can model the inaccuracy in the approximation error due to the finite sample size by including a prior additionally on the background mean and error covariance $p\left(\overline{\pmb{x}}_k^\mathrm{f}, \mathbf{B}_k^\mathrm{f}\right)$, as in the finite-size ensemble Kalman filter formalism of  \textcite{bocquet2011},  \textcite{bocquet2012} and \textcite{bocquet2015}.  Mathematical results demonstrate that covariance inflation can ameliorate the systematic underestimation of the variances due to model error in the presence of a low rank ensemble \parencite{grudzien2018}. In the presence of significant model error, the finite-size analysis is extended by the variant developed by \textcite{raanes2019adaptive}.

The linear transform of the ensemble matrix in Eq.~\eqref{eq:right_ensemble_update} is key to the efficiency of the (deterministic) EnKF presented above.  The ensemble-based cost function Hessian  $\boldsymbol{\Xi}_{\mathcal{J}_\mathrm{EnKF}} \in \mathbb{R}^{N_e \times N_e}$, where $N_e \ll N_x$ for typical geophysical models.  Particularly, the cost of computing the optimal weights as in Eq.~\eqref{eq:enkf_newton} and computing the transform matrix $\mathbf{T}$ in Eq.~\eqref{eq:right_transform} are both subordinate to the cost of the eigenvalue decomposition of the Hessian at $\mathcal{O}\left(N_e^3\right)$ floating point operations (flops), or to a randomised singular value decomposition \parencite{farchi2019}.  However, the extremely low ensemble size means that the correction to the forecast is restricted to the low-dimensional ensemble span, which may fail to correct directions of rapidly growing errors due to the rank deficiency.  While chaotic, dissipative dynamics implies that the background covariance $\mathbf{B}_k^{\mathrm{f/a}}$ has spectrum concentrated on a reduced rank subspace \parencite[and referenced therein]{carrassi2022}, covariance hybridisation \parencite{penny2017} or localisation \parencite{sakov2011} are used in practice to regularise the estimator's extreme rank deficiency, and the spurious empirical correlations that occur as a result of the degenerate sample size.

When extended to nonlinear dynamics, as in Eq.~\eqref{eq:state_model}, the EnKF can be seen to make an ensemble-based approximation to the EKF cost function, where the forecast ensemble is defined by
\begin{align}
\pmb{x}_{k,i}^{\mathrm{f}}:= \mathcal{M}_{k-1}\left(\pmb{x}_{k-1,i}^{\mathrm{a}} \right), & & \mathbf{E}_k^\mathrm{f}:= \begin{pmatrix}\pmb{x}_{k,1}^{\mathrm{f}} & \cdots & \pmb{x}_{k,N_e}^{\mathrm{f}} \end{pmatrix}.
\end{align}
The accuracy of the linear-Gaussian approximation to the dynamics of the nonlinear evolution of the ensemble, like the approximation of the EKF, depends strongly on the length of the forecast horizon $\Delta_t$. When the ensemble mean is a sufficiently accurate approximation of the mean state, and if the ensemble spread is of the same order as the error in the mean estimate \parencite{whitaker1998}, a similar approximation can be made for the ensemble evolution at first order as with Eq.~\eqref{eq:tangent_linear_model} but linearised about the ensemble mean as discussed later in Eq.~\eqref{eq:ensemble_linearization}.  Despite the similarity to the EKF, in the moderately nonlinear dynamics present in medium- to longer-range forecast horizons, the EnKF does not suffer from the same inaccuracy as the EKF in truncating the time-evolution at first order \parencite{evensen2003}.  The forecast ensemble members themselves evolve fully nonlinearly, tracking the higher-order dynamics, but the analysis update based on the Gaussian approximation becomes increasingly biased and fails to discriminate features like multi-modality of the posterior, even though this is occasionally an asset with higher-order statistical artefacts \parencite[see the discussion in][]{lawson2004}. 

\subsection{The Maximum Likelihood Ensemble Filter}
The EnKF filter analysis automatically accommodates weak nonlinearity in the state model without using the tangent linear model, as the filtering cost function has no underlying dependence on the state model. However, the EnKF analysis must be adjusted to account for nonlinearity in the observation model.   When a nonlinear observation operator is introduced, as in Eq.~\eqref{eq:observation_model}, the EnKF cost function can be re-written in the incremental analysis as
\begin{align}
\mathcal{J}_{\mathrm{EnKFI}}(\pmb{w}) := \frac{1}{2}(N_e - 1)\parallel \pmb{w}\parallel^2 +\frac{1}{2}\parallel \pmb{y}_k - \mathcal{H}_k\left(\hat{\pmb{x}}_k^\mathrm{i,f}\right)  - \mathbf{H}_k \mathbf{X}_k^\mathrm{f} \pmb{w}\parallel_{\mathbf{R}_k}^2,
\label{eq:incremental_ensemble_cost_weight_space}
\end{align}
where $\hat{\pmb{x}}_k^{i,\mathrm{f}}$ refers to the $i$-th iteration for the forecast mean.  Note that this does not refer to an estimate derived by an iterated ensemble forecast through the nonlinear state model -- rather, this is an iteration only with respect to the estimate of the optimal weights for the forecast perturbations.  If $\hat{\pmb{w}}^i$ is defined as the argmin of the cost function in Eq.~\eqref{eq:incremental_ensemble_cost_weight_space}, then the iterations of the ensemble mean are given as
\begin{align}
\hat{\pmb{x}}_k^{i+1,\mathrm{f}}:= \hat{\pmb{x}}^{i, \mathrm{f}}_k + \mathbf{X}_{k}^\mathrm{f} \hat{\pmb{w}}^{i} .
\end{align}
When $\parallel \hat{\pmb{w}}^i\parallel$ is sufficiently small, the optimisation terminates and the transform and the ensemble update can be performed as with the ETKF as in Eqs.~\eqref{eq:right_transform} and \eqref{eq:right_ensemble_update}. However, the direct, incremental approach above used the computation of the Jacobian of the observation operator $\mathbf{H}_k$.

The \emph{maximum likelihood ensemble filter} (MLEF) of \textcite{zupanski2005} and \textcite{zupanski2008} is an estimator designed to perform the above incremental analysis in the ETKF formalism, but without taking an explicit Taylor expansion for the observation operator above.  The method is designed to approximate the directional derivative of the nonlinear observation operator with respect to the ensemble perturbations,
\begin{align}\label{eq:3d_push_forward}
\mathbf{H}_k \mathbf{X}_k^\mathrm{f} := \nabla\vert_{\hat{\pmb{x}}^{i,\mathrm{f}}_{k}} \left[\mathcal{H}_k \right] \mathbf{X}^\mathrm{f}_k 
\end{align}
equivalent to computing the ensemble sensitivities of the map linearised about the ensemble-based mean.   This is often performed with an explicit finite-differences approximation between the ensemble members and the ensemble mean, mapped to the observation space.  Particularly, one may write
\begin{subequations}
\begin{align}
\mathcal{H}_k\left(\pmb{x}_k^{i, \mathrm{f}}\right)\approx \hat{\pmb{y}}^{i}_k &:=  \mathcal{H}_k \left(\hat{\pmb{x}}_k^{i,\mathrm{f}}\pmb{1}^\top   + \epsilon \mathbf{X}_k^\mathrm{f}\right) \pmb{1} / N_e\\
\mathbf{H}_k \mathbf{X}_k^\mathrm{f} \approx \widetilde{\mathbf{Y}}_k  &:= \frac{1}{\epsilon} \mathcal{H}_k \left(\hat{\pmb{x}}_k^{i,\mathrm{f}}\pmb{1}^\top   + \epsilon \mathbf{X}_k^\mathrm{f}\right) \left(\mathbf{I}_{N_e} - \pmb{1}\pmb{1}^\top / N_e \right),
\end{align}
\label{eq:finite_differences}
\end{subequations}
where $\epsilon$ is a small constant that re-scales the ensemble perturbations to approximate infinitesimals about the mean, and rescales the finite differences about the ensemble mean in the observation space.  This technique is used, e.g., in a modern form of the MLEF algorithm, based on the analysis of the iterative ensemble Kalman filter and smoother \parencite[see Section 6.7.2.1]{asch2016}.  The approximation above easily generalises to a 4D analysis, taking the directional derivative with respect to the state and observation models simultaneously, as in Eq.~\eqref{eq:4d_ensemble_finite_differences}.  However, in the smoothing problem, one can also extend the above analysis in terms of an alternating forwards filtering pass and backwards filtering pass to estimate the joint posterior.  This tutorial returns to the 4D analysis with the iterative ensemble Kalman filter and smoother in Section \ref{section:ienks} after an interlude on the retrospective analysis of the EnKS in the following section.

\subsection{The Ensemble Transform Kalman Smoother}
\label{section:enks}
The EnKS extends the filter analysis in the ETKF over the smoothing DAW by sequentially re-analysing past states with future observations with an additional filtering pass over the DAW backwards-in-time.  This analysis is performed retrospectively in the sense that the filter cycle of the ETKF is left unchanged, while an additional inner-loop of the DA cycle performs an update on the estimated lagged state ensembles within the DAW, stored in memory.  This can be formulated both for a fixed DAW and for fixed-lag smoothing, where only minor modifications are needed.  Consider here the algorithmically stationary DAW $\{t_1,\cdots ,t_L\}$ of fixed-lag smoothing, with a shift $S$ and lag $L$, and where it is assumed that  $S=1 \leq L$. The fixed-lag smoothing cycle of the EnKS begins by estimating the joint posterior pdf $p\left(\pmb{x}_{L:1}\vert \pmb{y}_{L:1}\right)$ recursively, given the joint posterior estimate over the last DAW $p\left(\pmb{x}_{L-1:0}\vert \pmb{y}_{L-1:0}\right)$. 

Given $p(\pmb{x}_{L-1:0}, \pmb{y}_{L-1:0})$, one can write the filtering pdf up to proportionality:
\begin{subequations}
\begin{align}
p(\pmb{x}_L \vert \pmb{y}_{L:0}) &\propto p(\pmb{y}_L \vert \pmb{x}_L, \pmb{y}_{L-1:0}) p(\pmb{x}_L , \pmb{y}_{L-1:0}) \\
&\propto \underbrace{p(\pmb{y}_L \vert \pmb{x}_L)}_{(i)} \underbrace{\int p(\pmb{x}_L \vert \pmb{x}_{L-1}) p(\pmb{x}_{L-1:0} \vert \pmb{y}_{L-1:0})\mathrm{d} \pmb{x}_{L-1:0}}_{(ii)},
\end{align}
\end{subequations}
as the product of (i) the likelihood of the observation $\pmb{y}_L$ given $\pmb{x}_L$; and (ii) the forecast for $\pmb{x}_L$ using the transition pdf on the last joint posterior estimate, marginalising out the past history of the model state $\pmb{x}_{L-1:0}$.  Recalling that $p(\pmb{x}_L \vert \pmb{y}_{L:1}) \propto p(\pmb{x}_L \vert \pmb{y}_{L:0})$, this provides a means to estimate the filter marginal of the joint posterior.  An alternating filtering pass, backwards-in-time, completes the smoothing cycle by estimating the joint posterior pdf $p(\pmb{x}_{L:1},\pmb{y}_{L:1})$.

Consider that the marginal smoother pdf is proportional to
\begin{subequations}
\begin{align}
p(\pmb{x}_{L-1} \vert \pmb{y}_{L:0}) &\propto p(\pmb{y}_L \vert \pmb{x}_{L-1}, \pmb{y}_{L-1:0})p(\pmb{x}_{L-1} , \pmb{y}_{L-1:0}) \\
&\propto \underbrace{p(\pmb{y}_L \vert \pmb{x}_{L-1})}_{(i)}\underbrace{p(\pmb{x}_{L-1} \vert\pmb{y}_{L-1:0})}_{(ii)},
\end{align}
\end{subequations}
where: (i) is the likelihood of the observation $\pmb{y}_L$ given the past state $\pmb{x}_{L-1}$; (ii) is the marginal pdf for $\pmb{x}_{L-1}$ from the last joint posterior.  The corresponding linear-Gaussian Bayesian MAP cost function is given for the retrospective analysis of the KS as
\begin{align}
\mathcal{J}_\mathrm{KS}(\pmb{x}_{L-1}) = \frac{1}{2} \parallel \pmb{x}_{L-1} - \overline{\pmb{x}}_{L-1|L-1}^\mathrm{s}\parallel^2_{\mathbf{B}_{L-1|L-1}^\mathrm{s}} +  \frac{1}{2} \parallel \pmb{y}_L - \mathbf{H}_L \mathbf{M}_L \pmb{x}_{L-1}\parallel_{\mathbf{R}_L}^2,
\end{align}
where $\overline{\pmb{x}}_{L-1|L-1}^\mathrm{s}$ and $\mathbf{B}^\mathrm{s}_{L-1|L-1}$ are the mean and covariance of the marginal smoother pdf $p(\pmb{x}_{L-1}\vert \pmb{y}_{L-1:0})$.  Define the matrix decomposition with the factorisation (e.g., a Cholesky decomposition):
\begin{align}
\mathbf{B}_{L-1|L-1}^\mathrm{s} = \boldsymbol{\Sigma}_{L-1|L-1}^\mathrm{s} \left(\boldsymbol{\Sigma}_{L-1|L-1}^\mathrm{s}\right)^\top,
\end{align}
and write $\pmb{x}_{L-1} = \overline{\pmb{x}}_{L-1|L-1}^\mathrm{s} + \boldsymbol{\Sigma}_{L-1|L-1}^\mathrm{s}\pmb{w}$, rendering the cost function as
\begin{subequations}
\begin{align}
\mathcal{J}_\mathrm{KS}(\pmb{w}) &= \frac{1}{2}\parallel \pmb{w}\parallel^2 + \frac{1}{2}\parallel \pmb{y}_L - \mathbf{H}_L \mathbf{M}_L (\overline{\pmb{x}}_{L-1|L-1}^\mathrm{s}  +  \boldsymbol{\Sigma}_{L-1|L-1}^\mathrm{s} \pmb{w})\parallel_{\mathbf{R}_{L}}^2\\
&=\frac{1}{2} \parallel\pmb{w}\parallel^2 + \frac{1}{2}\parallel \pmb{y}_L - \mathbf{H}_L \overline{\pmb{x}}_{L}^\mathrm{f} - \mathbf{H}_L \boldsymbol{\Sigma}_{L}^\mathrm{f}\pmb{w} \parallel_{\mathbf{R}_L}^2.
\end{align}
\label{eq:retrospective_exact}
\end{subequations}
Let $\overline{\pmb{w}}$ now denote the argmin of Eq.~\eqref{eq:retrospective_exact}.  It is important to recognise that
\begin{subequations}
\begin{align}
\pmb{x}_L := & \mathbf{M}_L \left( \overline{\pmb{x}}_{L-1|L-1}^\mathrm{s} + \boldsymbol{\Sigma}_{L-1|L-1}^\mathrm{s} \pmb{w}\right) \\
= &\overline{\pmb{x}}_L^\mathrm{f} + \boldsymbol{\Sigma}_L^\mathrm{f} \pmb{w},
\end{align}
\end{subequations}
such that the argmin for the smoothing problem $\overline{\pmb{w}}$ is also the argmin for the filtering MAP analysis.  

The ensemble-based approximation,
\begin{subequations}
\begin{align}
\pmb{x}_{L-1}&= \hat{\pmb{x}}_{L-1|L-1}^\mathrm{s} + \mathbf{X}_{L-1|L-1}^\mathrm{s}\pmb{w},\\
\mathcal{J}_\mathrm{EnKS}(\pmb{w}) &=\frac{1}{2}(N_e - 1) \parallel\pmb{w}\parallel^2 + \frac{1}{2}\parallel \pmb{y}_L - \mathbf{H}_L\hat{\pmb{x}}^\mathrm{f}_L - \mathbf{H}_L\mathbf{X}_L^\mathrm{f}\pmb{w} \parallel^2,
\end{align}
\label{eq:retrospective_ensemble}
\end{subequations}
to the exact smoother cost function in Eq.~\eqref{eq:retrospective_exact} yields the retrospective analysis of the EnKS as
\begin{subequations}
\label{eq:EnKS_analysis}
\begin{align}
\hat{\pmb{w}} &:= \pmb{0} - \boldsymbol{\Xi}_{\mathcal{J}_{\mathrm{EnKS}}}^{-1} \nabla \mathcal{J}_{\mathrm{EnKS}}|_{\pmb{w}=\pmb{0}},\\
\mathbf{T} &: = \boldsymbol{\Xi}_{\mathcal{J}_{\mathrm{EnKS}}}^{-\frac{1}{2}},\\
\begin{split}
\mathbf{E}_{L-1|L}^\mathrm{s} &= \hat{\pmb{x}}_{L-1| L-1}^\mathrm{s} \pmb{1}^\top + \mathbf{X}_{L-1| L-1}^\mathrm{s}\left(\hat{\pmb{w}} \pmb{1}^\top + \sqrt{N_e - 1} \mathbf{T}\mathbf{U}\right), \\
& = \mathbf{E}_{L-1|L-1}^\mathrm{s} \boldsymbol{\Psi}_L.
\end{split}\label{eq:enks_retrospective_conditioning}
\end{align}
\end{subequations}
where $\boldsymbol{\Psi}_L$ is the ensemble transform as defined for the filtering update as in Eq.~\eqref{eq:right_ensemble_update}.

The above equations generalise for arbitrary indices $k|L$ over the DAW, providing the complete description of the inner-loop between each filter cycle of the EnKS. After each new observation is assimilated with the ETKF analysis step, a smoother inner-loop makes a backwards pass over the DAW applying the transform and the weights of the ETKF filter update to each past ensemble state stored in memory. This analysis easily generalises to the case where there is a shift of the DAW with $S>1$, though the EnKS alternating forwards and backwards filtering passes must be performed in sequence over the observations, ordered-in-time, rather than making a global analysis over $\pmb{y}_{L:L-S+1}$.  Finally, this easily extends to accommodate a nonlinear observation model by using the MLEF filtering step to obtain the optimal ensemble transform, and by applying this recursively backwards-in-time to the lagged ensemble states. 

\begin{figure*}
\center
\includegraphics[width=\linewidth]{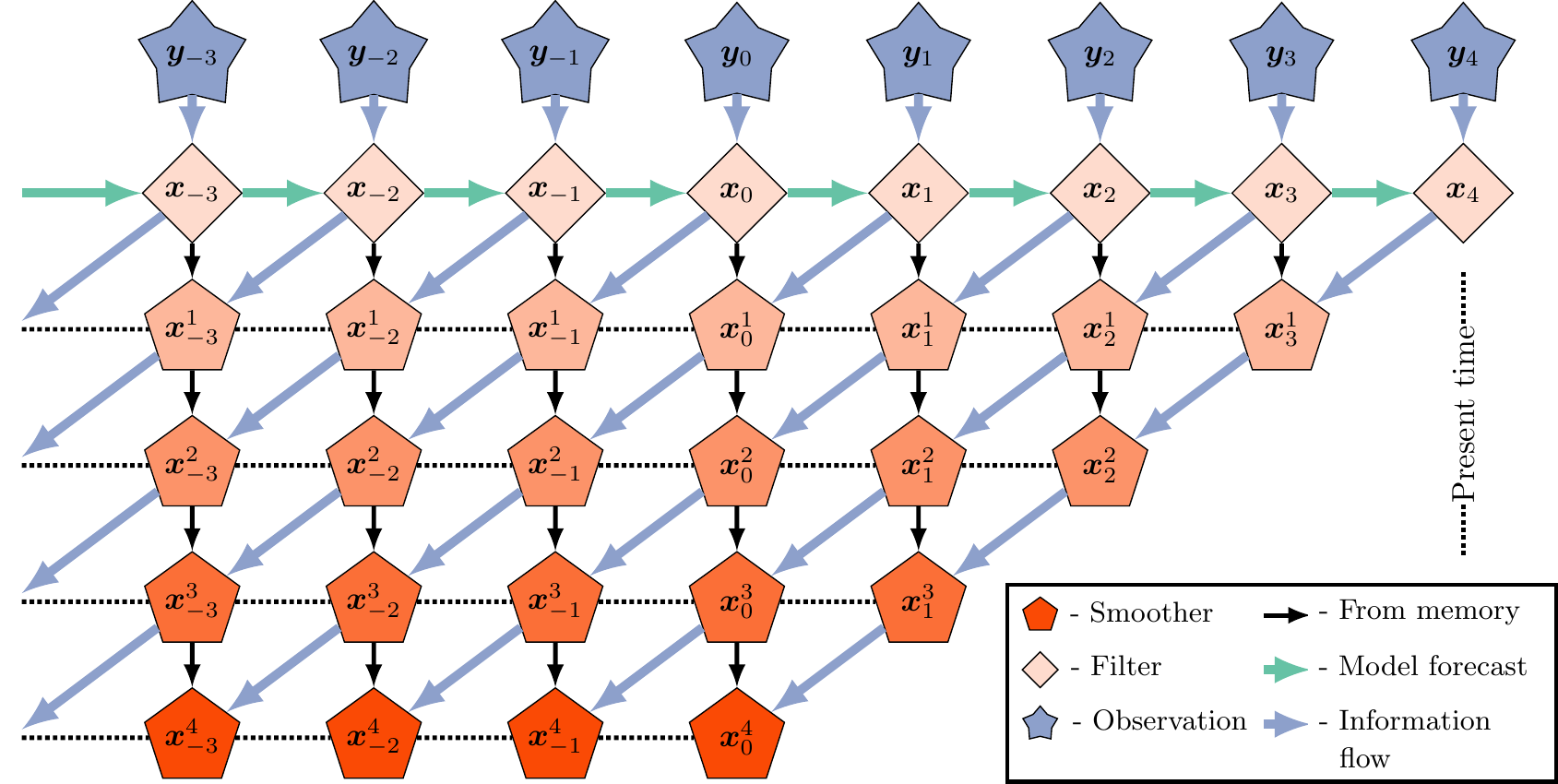}
\caption{$L=4$ (lag), $S=1$ (shift) EnKS.  Observations are assimilated sequentially via the filter cost function and a retrospective re-analysis is applied to all ensemble states within the lag window stored in memory.  Source: \textcite{grudzien2021fast}, adapted from \textcite{asch2016}.}
\label{fig:enks_diagram}
\end{figure*}

A schematic of the EnKS cycle for a lag of $L=4$ and a shift of $S=1$ is pictured in Fig.~\ref{fig:enks_diagram}.  Time moves forwards from left to right in the horizontal axis with a step size of $\Delta_t$.  At each analysis time, the ensemble forecast from the last filter pdf is combined with the observation to produce the ensemble transform update.  This transform is then utilised to produce the posterior estimate for all lagged ensemble states, conditioned on the new observation.  The information in the posterior estimate thus flows in reverse time to the lagged states stored in memory, but the information flow is unidirectional in this scheme.  This type of retrospective analysis maintains the computational cost of the fixed-lag EnKS described above at a comparable level to the EnKF, with the only significant additional cost being the storage of the ensemble at lagged times to be re-analysed.  

\subsection{The Iterative Ensemble Kalman Filter and Smoother}
\label{section:ienks}
While the EnKS is computationally efficient, its retrospective analysis has some drawbacks compared to the 4D analysis in terms of forecast accuracy.  The 4D fixed-lag smoothing analysis re-initialises each cycle with a re-analysed estimate for the initial data, transmitting the observations' information forwards-in-time through the nonlinear dynamics.  The challenge with the 4D analysis in the absence of the tangent linear and adjoint models is in devising how to efficiently and accurately compute the gradient of the 4D cost function.  Building on \textcite{zupanski2005,liu2008}, the analysis of the \emph{iterative ensemble Kalman filter} (IEnKF) and the \emph{iterative ensemble Kalman smoother} (IEnKS) extends the ensemble transform method of the ETKF to iteratively optimise the 4D cost function.

Recall the quadratic cost function at the basis of incremental 4D-VAR, Eq.~\eqref{eq:incremental_4dvar} -- making an ensemble-based approximation for the background mean and error covariance, and writing the model state as an ensemble perturbation of the $i$-th proposal for the ensemble mean, one has
\begin{subequations}
\begin{align}
&\begin{split}
\mathcal{J}_{\mathrm{IEnKS}}(\pmb{w}):=& \frac{1}{2} \parallel \hat{\pmb{x}}^i_0 - \hat{\pmb{x}}^i_0 -\mathbf{X}_0 \pmb{w} \parallel_{\mathbf{P}_0}^2\\
& + \frac{1}{2} \sum_{k=1}^L \parallel \pmb{y}_k - \mathcal{H}_k \circ \mathcal{M}_{k:1} \left( \hat{\pmb{x}}^i_0\right) - \mathbf{H}_k\mathbf{M}_{k:1}\mathbf{X}_0 \pmb{w} \parallel_{\mathbf{R}_k}^2 
\end{split}\\
 & 
\begin{split}
=& \frac{1}{2} \left( N_e - 1\right) \parallel\pmb{w} \parallel^2\\
& + \frac{1}{2} \sum_{k=1}^L \parallel \pmb{y}_k - \mathcal{H}_k \circ \mathcal{M}_{k:1} \left( \hat{\pmb{x}}^i_0\right) - \mathbf{H}_k\mathbf{M}_{k:1}\mathbf{X}_0 \pmb{w} \parallel_{\mathbf{R}_k}^2 .
\end{split}
\end{align}
\end{subequations}
Applying the finite differences approximation as given in Eq.~\eqref{eq:finite_differences}, but with respect to the composition of the nonlinear observation operator with the nonlinear state model
\begin{subequations}
\begin{align}
\mathcal{H}_k \circ \mathcal{M}_{k:1} \left(\hat{\pmb{x}}_0^i\right) \approx \hat{\pmb{y}}_k&:=  \mathcal{H}_k \circ \mathcal{M}_{k:1}\left(\hat{\pmb{x}}^{i}_0 + \mathbf{X}_0\epsilon \right)\pmb{1}/N_e ,\\
\mathbf{H}_k \mathbf{M}_{k:1}\mathbf{X}_0 \approx \widetilde{\mathbf{Y}}_k &:=\frac{1}{\epsilon}\mathcal{H}_k \circ \mathcal{M}_{k:1}\left(\hat{\pmb{x}}^{i}_0 + \mathbf{X}_0\epsilon\right)\left(\mathbf{I}_{N_e} - \pmb{1}\pmb{1}^\top/N_e \right),
\end{align}
\label{eq:4d_ensemble_finite_differences}
\end{subequations}
this sketches the `bundle' formulation of the IEnKS \parencite{bocquet2013,bocquet2014}.  The above indexing refers to the case where the smoothing problem is performed offline, with a fixed DAW, though a similar indexing to Eq.~\eqref{eq:4d_map_cost_sequential} gives the sequential form over newly incoming observations.  

The sequential form of the IEnKS can be treated as a nonlinear sequential filter, which is the purpose that the method was originally devised for.  Indeed, setting the number of lagged states $L=1$, this provides a direct extension of the EnKF cost function, Eq.~\eqref{eq:ensemble_cost_weight_space}, but where there is an additional dependence on the initial conditions for the ensemble forecast.  This lag-1 iterative filtering scheme is called the IEnKF \parencite{sakov2012,bocquet2012}, which formed the original basis for the IEnKS.  Modern forms of the IEnKF/S analysis furthermore include the treatment of model errors, like weak-constraint 4D-VAR \parencite{sakov2018a,sakov2018b,fillion2020}.  Alternative formulations of this analysis, based on the original stochastic, perturbed observation EnKF \parencite{evensen1994,burgers1998}, also have a parallel development in the \emph{ensemble randomised maximum likelihood method} (EnRML) of \textcite{gu2007,chen2012,raanes2019revising}.  A similar ensemble-variational estimator based on the EnKF analysis is the \emph{ensemble Kalman inversion} (EKI) of \textcite{iglesias2013, schillings2018, kovachki2019}.

\subsection{The Single-Iteration Ensemble Kalman Smoother}
The sequential smoothing analysis so far has presented two classical approaches: (i) a 4D approach using the global analysis of all new observations available within a DAW at once, optimising an initial condition for the lagged model state; and (ii) the 3D approach based upon alternating forwards and backwards filtering passes over the DAW.  Each of these approaches has strengths and weaknesses in a computational cost / forecast accuracy trade off.  Particularly, the 4D approach, as in the IEnKS, benefits from the improved estimate of the initial condition when producing the subsequent forecast statistics in a shifted DAW; however, each step of the iterative optimisation comes at the cost of simulating the ensemble forecast over the entire DAW in the fully nonlinear state model, which is typically the greatest numerical expense in geophysical DA.  On the other hand, the 3D approach of the classical EnKS benefits from a low computational cost, requiring only a single ensemble simulation of the nonlinear forecast model over the DAW, while the retrospective analysis of lagged states is performed with the filtering transform without requiring additional model simulations (though at a potentially large memory storage cost).  It should be noted that in the perfect, linear-Gaussian model Bayesian analysis, both approaches produce equivalent estimates of the joint posterior.  However, when nonlinearity is present in the DA cycle the approaches produce distinct estimates -- for this reason, the source of the nonlinearity in the DA cycle is of an important practical concern.

Consider the situation in which the forecast horizon $\Delta_t$ is short, so that the ensemble time-evolution is weakly nonlinear, and where the ensemble mean is a good approximation for the mean of the filtered distribution.  In this case, the model forecast dynamics are well-approximated by the tangent linear evolution of a perturbation about the ensemble mean, i.e.,
\begin{align}
\pmb{x}_k = \mathcal{M}_k\left(\pmb{x}_{k-1}\right) \approx \mathcal{M}_k \left(\hat{\pmb{x}}_{k-1}\right) + \mathbf{M}_k \pmb{\delta}_{k-1} .\label{eq:ensemble_linearization}
\end{align}
For such short-range forecasts, nonlinearity in the observation-analysis-forecast cycle may instead be dominated by the nonlinearity in the observation operator $\mathcal{H}_k$, or in the optimisation of hyper-parameters of the filtering cost function, and not by the model forecast itself.  In this situation, an iterative simulation of the forecast dynamics as in the 4D approach may not produce a cost-effective reduction in the analysis error as compared to, e.g., the MLEF filter analysis optimising the filtering cost function alone.  However, one may still obtain the benefits of re-initialisation of the forecast with a re-analysed prior by using a simple hybridisation of the 3D and 4D smoothing analyses.  

Specifically, in a given fixed-lag smoothing cycle, one may iteratively optimise the sequential filtering cost functions for a given DAW corresponding to the new observations at times $\{t_{L-S+1},\cdots , t_{L}\}$ as with the MLEF.  These filtering ensemble transforms not only condition the ensemble at the corresponding observation time instance, but also produce a retrospective analysis of a lagged state as in Eq.~\eqref{eq:enks_retrospective_conditioning} with the EnKS.   However, when the DAW itself is shifted, one does not need to produce a forecast from the latest filtered ensemble -- instead one can initialise the next ensemble forecast with the lagged, retrospectively smoothed ensemble at the beginning of the last DAW.  This hybrid approach utilising the retrospective analysis, as in the classical EnKS, and the ensemble simulation over the lagged states while shifting the DAW, as in the 4D analysis of the IEnKS, was recently developed by \textcite{grudzien2021fast}, and is called the \emph{single-iteration ensemble Kalman smoother} (SIEnKS).  The SIEnKS is named as such because it produces its forecast, filter and re-analysed smoother statistics with a single iteration of the ensemble simulation over the DAW in a fully consistent Bayesian analysis. By doing so, it seeks to minimise the leading order cost of EnVAR smoothing, i.e., the ensemble simulation in the nonlinear forecast model.  However, the estimator is free to iteratively optimise the filter cost function for any single observation vector without additional iterations of the ensemble simulation.  In observation-analysis-forecast cycles in which the forecast error dynamics are weakly nonlinear, yet other aspects of the cycle are moderately to strongly nonlinear, this scheme is shown to produce a forecast accuracy comparable to, and at times better than, the 4D approach but with an overall lower leading-order computational burden \parencite{grudzien2021fast}.

\begin{figure*}
\center
\includegraphics[width=\linewidth]{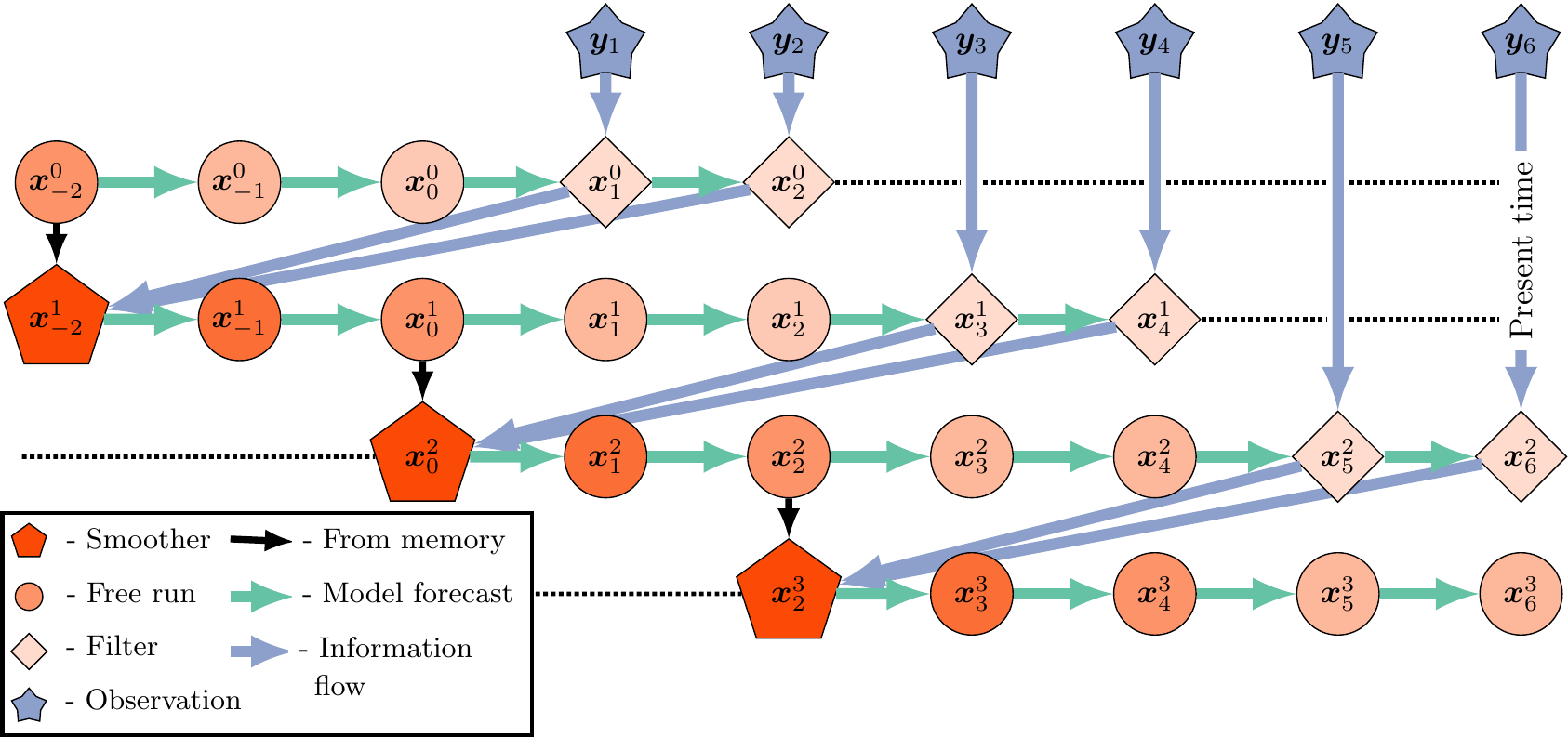}
\caption{$L=4$ (lag), $S=2$ (shift) SIEnKS diagram.  An initial condition from the last smoothing cycle initialises a forecast simulation over the current DAW of $L=4$ states.  New observations entering the DAW are assimilated sequentially via the filter cost function.  After each filter analysis, a retrospective re-analysis is applied to the initial ensemble, and this re-analysed initial condition is evolved via the model $S$ analysis times forwards to begin the next cycle. Source: \textcite{grudzien2021fast}.}
\label{fig:sienks_diagram}
\end{figure*}

A schematic of the SIEnKS cycle for a lag of $L=4$ and a shift of $S=2$ is pictured in Fig. \ref{fig:sienks_diagram}.  This demonstrates how the sequential analysis of the filter cost function and sequential, retrospective re-analysis for each incoming observation differs from the global analysis of the 4D approach of the IEnKS. Other well-known DA schemes combining a retrospective re-analysis and re-initialisation of the ensemble forecast include the \emph{Running In Place} (RIP) smoother of  \textcite{kalnay2010, yang2013} and the \emph{One Step Ahead} (OSA) smoother of  \textcite{desbouvries2011} and \textcite{ait2022}.  
It can be shown that, with an ETKF style filter analysis, a single iteration of the ensemble over the DAW, a perfect model assumption and a lag of $L=S=1$, the SIEnKS, RIP and OSA smoothers all coincide \parencite{grudzien2021fast}.

\section{Machine Learning and Data Assimilation in a Bayesian Perspective}

This final section demonstrates how the Bayesian DA framework can be extended to estimate more than the state vector, including key parameters of the model and the observation-analysis-forecast cycle.  In particular, the Bayesian framework can be used to formulate techniques for learning both the state vector and part of, if not the full, dynamical model.  If this objective was always in the scope of classical DA, it was actually made possible, significantly beyond linear regression, by the introduction of \emph{Machine Learning} (ML) techniques to supplement traditional DA schemes.

Incorporating ML into DA algorithms was suggested quite early by \textcite{hsieh1998}, who clearly advocated the use of \emph{Neural Networks} (NNs) in a variational DA framework.
More recently this was put forward and illustrated by \textcite{abarbanel2018} and \textcite{bocquet2019} with a derivation of the generalised cost function from Bayes' law, 
showing that the classical ML cost function for learning a surrogate model of the dynamical model is a limiting case of the DA framework, as in Eq.~\eqref{eq:ml_cost}.
This formalism was furthermore generalised by \textcite{bocquet2020}, using classical DA notation, for learning the state vector, the dynamical model and the error statistics attached to this full retrieval.  This section follows this latest paper, showing how to derive the cost function of the generalised estimation problem.   Note that going even beyond this approach, attempting to learn the optimisation scheme of the cost function, or even the full DA procedure, an approach called \emph{end-to-end} in ML, is a subject of active investigations \parencite{fablet2021,peyron2021}.

\subsection{Prior Error Statistics}

For the sake of simplicity, again assume Gaussian statistics for the observation errors, $p(\pmb{y}_k | \pmb{x}_k ) = n(\pmb{y}_k | \pmb{x}_k, \mathbf{R}_k)$,
where the observation error covariance matrices $\bR_{L:1}:=\left\{\bR_L, \bR_{L-1}, \cdots, \bR_1\right\}$ are supposed to be known.
The dynamical model is meant to be learned or approximated and thus stands as a surrogate model for the unknown true physical dynamics. Assuming that the model does not explicitly depend on time, its resolvent is defined by
\begin{align}
\pmb{x}_k := \bF^{k}_\bA \left(\pmb{x}_{k-1}\right) + \pmb{\eta}_k,
\label{eq:surrogate_model}
\end{align}
depending on a (possibly very large) set of parameters $\bA$.  Prototypically, $\mathbf{A}$ represents the weights and biases of a NN, which are learned from the observations alongside the state vectors within the DAW.
The distribution for model error, i.e., $\boldsymbol{\eta}_k$ in Eq.~\eqref{eq:surrogate_model}, is also assumed Gaussian such that 
\be
p(\pmb{x}_k | \pmb{x}_{k-1}, \mathbf{A}, \mathbf{Q}_k) := n\left(\pmb{x}_k | \mathbf{F}_\mathbf{A}^k(\pmb{x}_{k-1}), \mathbf{Q}_k\right),
\ee
where $\bQ_{L:1}:=\left\{\bQ_L, \bQ_{L-1}, \cdots, \bQ_1\right\}$ are not necessarily known.  Further assume that these Gaussian errors are white-in-time and that the observation and model errors are mutually independent.

Note the intriguing status of $\bQ_{L:1}$ since it depends, a posteriori, on how well the surrogate model is estimated. This calls for an adaptive estimation of the model error statistics $\bQ_{L:1}$ as the surrogate model, parametrised by $\mathbf{A}$, is better approximated.

\subsection{Joint Estimation of the Model, its Error Statistics, and the State Trajectory}

Following the Bayesian formalism, one may form a MAP estimate for the joint pdf in $\bA$ and $\pmb{x}_{L:0}$ conditioned on the observations, together with the model error statistics.
This generalised conditional pdf is expressed in the hierarchy
\begin{align}
  p(\bA, \bQ_{L:1}, \pmb{x}_{L:0} | &\pmb{y}_{k:1}, \bR_{K:0}) = \nn
&  \frac{p(\pmb{y}_{L:1}|\pmb{x}_{L:0},\bR_{L:0} )p(\pmb{x}_{L:0}|\bA, \bQ_{L:1})p(\bA, \bQ_{L:1})}{p(\pmb{y}_{L:1}, \bR_{L:0})} ,
\end{align}
where the mutual independence of the observation and model error is used. Once again, this remarkably stresses how powerful and general the Bayesian framework can be.

The first term in the numerator of the right-hand-side is the usual likelihood of the observations.
The second term in the numerator is the prior on the trajectory, given a known model and known model error statistics.  The final term of the numerator is the joint prior of the model and the model error statistics, described as a \emph{hyperprior}. The associated cost function is derived proportional to
\begin{align}
  \mathcal{J}_{\mathrm{DA}\text{-}\mathrm{ML}}(\bA,\pmb{x}_{L:0},\bQ_{L:1}) = & -\log \, p(\bA, \bQ_{L:1}, \pmb{x}_{L:0} | \pmb{y}_{L:1}, \bR_{L:1}) \nn
  =& \frac{1}{2}\sum_{k=1}^{L} \left\{ \left\| \pmb{y}_k-\mathcal{H}_k(\pmb{x}_k) \right\|^2_{\mathbf{R}_k} + \log \left(\left|\bR_k \right|\right) \right\} \nn
  &  + \frac{1}{2}\sum_{k=1}^{L} \left\{ \left\| \pmb{x}_k-\mathbf{F}^{k}_\bA(\pmb{x}_{k-1}) \right\|^2_{\mathbf{Q}_k}
+ \log \left(\left|\bQ_k \right|\right) \right\} \nn
& - \log \, p(\pmb{x}_0, \bA, \bQ_{L:1}) .
\label{eq:da_ml_cost}
\end{align}
Note the resemblance of \eqref{eq:da_ml_cost} with the weak-constraint 4D-VAR cost function of classical DA \parencite{tremolet2006}.  Very importantly, this Bayesian formulation allows for a rigorous treatment of partial and noisy observations.
The classical ML cost function that uses noiseless, complete observations of the physical system is derived from Eq.~\eqref{eq:da_ml_cost},
assuming that $\bQ_k$ is known, $\bH_k \equiv \mathbf{I}_{N_x}$, and setting $\bR_k$ to go to $\bzero$. Associating the initial data as $\pmb{y}_0 \longleftarrow \pmb{x}_0$, $\mathcal{J}_{\mathrm{DA}\text{-}\mathrm{ML}}(\bA,\pmb{x}_{L:0},\bQ_{L:1})$ becomes in its limit
\begin{align}
  \mathcal{J}_\mathrm{ML}(\bA) = \frac{1}{2}\sum_{k=1}^{L} \left\| \pmb{y}_k-\mathbf{F}^{k}_\bA(\pmb{y}_{k-1}) \right\|^2_{\mathbf{Q}_k}  - \log p(\pmb{y}_0,\bA,\bQ_{L:1}).
  \label{eq:ml_cost}
\end{align}
Such connections between ML and DA were first highlighted by \textcite{hsieh1998,abarbanel2018,bocquet2019}.

\begin{figure}[ht]
\includegraphics[width=\linewidth]{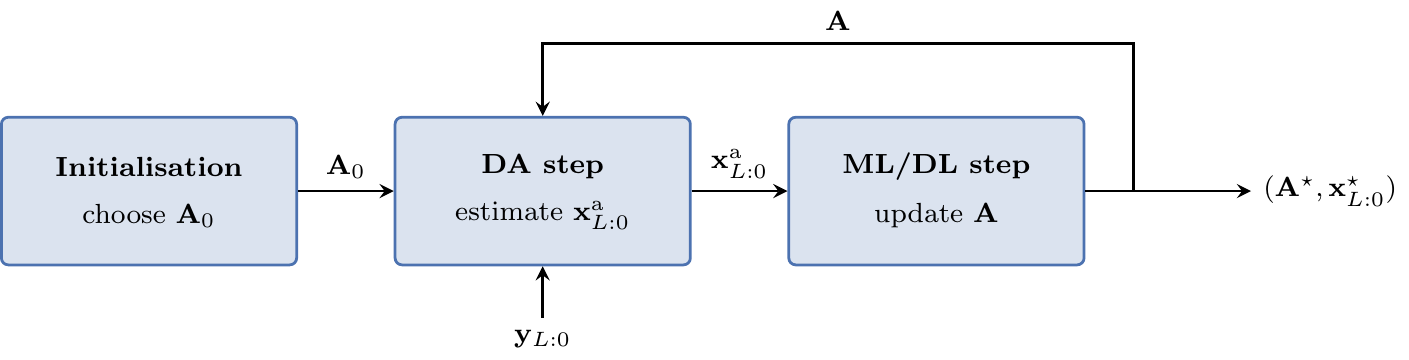}
\caption{\label{fig:loop} Estimation of both model and state trajectory using coordinate descent by alternately optimising on the state trajectory using DA and on the NN model parameters using ML/DL. The iterative loop stops when an accuracy criterion is met. Source: \textcite{farchi2021machine}}
\end{figure}

Solving the combined DA/ML problem, i.e. minimising Eq.~\eqref{eq:da_ml_cost} leads to several key remarks.  First, as mentioned earlier, this formalism allows one to learn a surrogate model of the true dynamics using partial, possibly sparse and noisy observations, as opposed to off-the-shelf ML techniques. This is obviously critical for geophysical systems.  Second, in this framework, one fundamentally looks for a stochastically additive surrogate model of the form in Eq.~\eqref{eq:surrogate_model} rather than a deterministic surrogate model, since $\ceta_{k}$ is drawn from the normal distribution of covariance matrices $\bQ_{k}$.  Third, there are many ways to carry out this minimisation, as discussed by \textcite{bocquet2019,bocquet2020}. Because the state vectors and the model parameters $\bA$ are of fundamentally different nature, and yet are statistically interdependent, one idea is to minimise Eq.~\eqref{eq:da_ml_cost} through a coordinate descent, i.e. alternating minimisations on $\bA$ and $\pmb{x}_{L:0}$, as illustrated by Fig.~\ref{fig:loop}. This was first suggested and successfully implemented by \textcite{brajard2020}, by using an EnKF for the assimilation step (i.e. minimising on $\pmb{x}_{L:0}$), and a \emph{Deep Learning} (DL) optimiser for the ML step (i.e., minimising on $\bA$).  This work was extended to 4D variational analysis using a 4D-VAR assimilation step with a DL optimiser for the ML step by \textcite{farchi2021machine}.

\subsection{Joint Estimation of the Model and the Error Statistics}

A slightly different objective is to obtain a MAP estimate for the surrogate model, irrespective of any model state realisation, i.e., one is interested in the MAP of the marginal conditional pdf
\be
\label{eq:p2}
p(\bA, \bQ_{L:1} | \pmb{y}_{L:1}, \bR_{L:1}) = \int \!  p(\bA, \bQ_{L:1}, \pmb{x}_{L:0} | \pmb{y}_{L:1}, \bR_{L:1}) \, \mathrm{d}\pmb{x}_{L:0} ,
\ee
which is theoretically obtained by minimising
\be
\label{eq:J2}
\mathcal{J}(\bA, \bQ_{L:1}) = -\log \, p(\bA, \bQ_{L:1} | \pmb{y}_{L:1}, \bR_{L:1}).
\ee
As pointed out by \textcite{bocquet2019}, the marginal pdf \eqref{eq:p2} can be approximately related to the joint pdf through a Laplace approximation of the integral.
Here, however, one is interested in the full solution to this problem. This can be solved numerically by using the \emph{Expectation-Maximisation} (EM) statistical algorithm \parencite{dempster1977} jointly with the variational numerical solution of Eq.~\eqref{eq:da_ml_cost}. This was suggested in \textcite{ghahramani1999,nguyen2019} and implemented and validated by \textcite{bocquet2020}.

\section{Conclusions}

Unifying techniques from statistical estimation, nonlinear optimisation and even machine learning, the Bayesian approach to DA  provides a consistent treatment of a variety of topics in the classical DA problem, discussed throughout this tutorial.  Furthermore, the Bayesian approach can be reasonably extended to treat additional challenges not fully considered in this tutorial, such as estimating process model parameters, handling significant modelling errors, and including the optimisation of various hyper-parameters of the observation-analysis-forecast cycle.  Estimation may include even learning the dynamical process model itself, as in the emerging topic of DA/ML hybrid algorithms. In this context especially, a Bayesian analysis provides a coherent treatment of the problem where either DA or ML techniques themselves may be insufficient.  This is particularly relevant where surrogate ML models are used to augment traditional physics-based dynamical models for, e.g., simulating unresolved dynamics at scales too fine to be dynamically represented.  The hybrid approach in the Bayesian analysis provides a means to combine the dynamical and surrogate simulations with real-world observations, and to produce an analysis of the state and parameters for which subsequent simulations depend. This tutorial thus presents one framework for interpreting the DA problem, both in its classical formulation and in the directions of the current state-of-the-art.  In surveying a variety of widely used DA schemes, the key message of this tutorial is how the Bayesian analysis provides a consistent framework for the estimation problem and how this allows one to formulate its solution in a variety of ways to exploit the operational challenges in the geosciences.  

\section*{Acknowledgements}
This tutorial benefited from the opportunity to lecture on data assimilation at: (i) the CIMPA Research School on \emph{Mathematics of Climate Science} in Kigali, Rawanda (28 June - 9 July, 2021); (ii) Utrecht University's Research school on
\emph{Data science and beyond: data assimilation with elements of machine learning} in Utrecht, Netherlands (23 August - 27 August, 2021); (iii) the Joint ICTP-IUGG Workshop on \emph{Data Assimilation and Inverse Problems in Geophysical Sciences} in Trieste, Italy (18 October - 29 October 2021); and (iv) at the University of Nevada, Reno, in Fall Semester 2021, where the combined lecture notes formed the basis of this manuscript.  The authors would like to thank everyone, but especially the students, the named reviewer Polly Smith, and the editor Alik Ismail-Zadeh, who supported the development of these notes with their many and highly relevant suggestions. Finally, this tutorial also benefited from many fruitful discussions with Alberto Carrassi, Alban Farchi, Patrick Raanes, Julien Brajard, and Laurent Bertino. CEREA is a member of Institut Pierre-Simon Laplace (IPSL).

\printbibliography

\end{document}

%% file: macro.tex
\newcommand\be{\begin{equation}}
\newcommand\ee{\end{equation}}
\newcommand\nn{\nonumber\\}

\newcommand\bA{\mathbf{A}}

\newcommand\bF{\mathbf{F}}
\newcommand\bQ{\mathbf{Q}}

\newcommand\bH{\mathbf{H}}
\newcommand\bR{\mathbf{R}}

\newcommand\bzero{\mathbf{0}}

\newcommand\ceta{\boldsymbol \eta}